\documentclass[twoside]{article}
\usepackage{amssymb}
\catcode`\@=11

\long\def\@makefntext#1{
\protect\noindent \hbox to 3.2pt {\hskip-.9pt  
$^{{\eightrm\@thefnmark}}$\hfil}#1\hfill}		

\def\@makefnmark{\hbox to 0pt{$^{\@thefnmark}$\hss}}	
	
\def\ps@myheadings{\let\@mkboth\@gobbletwo
\def\@oddhead{\hbox{}
\rightmark\hfil\eightrm\thepage}   
\def\@oddfoot{}\def\@evenhead{\eightrm\thepage\hfil
\leftmark\hbox{}}\def\@evenfoot{}
\def\sectionmark##1{}\def\subsectionmark##1{}}

\renewcommand{\section}[1] {\vspace{12pt}\addtocounter{section}{1}
\setcounter{Theorem}{0}\setcounter{equation}{0}
\setcounter{subsection}{0}\noindent
        {\tenbf\thesection. #1}\par\vspace{5pt}}
\renewcommand{\subsection}[1] {\vspace{12pt}\addtocounter{subsection}{1}
        \noindent
        {\bf\thesubsection. {\kern1pt \bfit #1}}\par\vspace{5pt}}
\newcommand{\nonumsection}[1] {\vspace{12pt}\noindent{\tenbf #1}
        \par\vspace{5pt}}

\newtheorem{Theorem}{Theorem}[section] 
\newtheorem{Lemma}[Theorem]{Lemma}
\newtheorem{Remark}[Theorem]{Remark}
\newtheorem{Corollary}[Theorem]{Corollary}
\newtheorem{Proposition}[Theorem]{Proposition}
\newtheorem{Example}[Theorem]{Example}

\newcommand{\textlineskip}{\baselineskip=13pt}
\newcommand{\smalllineskip}{\baselineskip=10pt}

\newcommand{\copyrightheading}[1]
     {\vspace*{-2.5cm}\smalllineskip{\flushleft
     {\footnotesize Communications in Contemporary Mathematics, #1}\\
     {\footnotesize \copyright\kern2pt World Scientific Publishing
         Company}\\
         }}

\newcommand{\publisher}[2]{{\begin{center}\footnotesize\smalllineskip
        Received #1\\
        Revised #2
        \end{center}
        }}

\def\abstracts#1#2#3{{
        \centering{\begin{minipage}{4.5in}\footnotesize\baselineskip=10pt
        \parindent=0pt #1\par
        \parindent=15pt #2\par
        \parindent=15pt #3
        \end{minipage}}\par}}

\newcommand{\bibit}{\nineit}
\newcommand{\bibbf}{\ninebf}
\renewenvironment{thebibliography}[1]
        {\frenchspacing
         \ninerm\baselineskip=11pt
         \begin{list}{\arabic{enumi}.}
        {\usecounter{enumi}\setlength{\parsep}{0pt}
        \setlength{\leftmargin 17pt}{\rightmargin 0pt}   
         \setlength{\itemsep}{0pt} \settowidth
        {\labelwidth}{#1.}\sloppy}}{\end{list}}

\def\fpage#1{\begingroup
\voffset=.3in
\thispagestyle{empty}\begin{table}[b]\centerline{\footnotesize #1}
        \end{table}\endgroup}

\def\runninghead#1#2{\pagestyle{myheadings}
\markboth{{\protect\footnotesize\it{\quad #1}}\hfill}
{\hfill{\protect\footnotesize\it{#2\quad}}}}
\font\tenbf=cmbx10
\font\bfit=cmbxti10 at 10pt
\font\ninerm=cmr9
\font\nineit=cmti9
\font\ninebf=cmbx9
\font\eightrm=cmr8

\def\qed{\hbox{${\vcenter{\vbox{                        
   \hrule height 0.4pt\hbox{\vrule width 0.4pt height 6pt
   \kern5pt\vrule width 0.4pt}\hrule height 0.4pt}}}$}\vspace{7pt}}
\def\R{{\mathbb R}}
\def\N{{\mathbb N}}
\def\S{{\mathbb S}}
\def\Spheres{{\mathcal S}}
\def\E{{\mathcal E}}
\def\B{{\mathcal B}}
\def\D{{\mathcal D}}
\def\V{{\mathcal V}}
\def\disc{D}

\def\div{{\rm{div}}}

\def\Proof{{\noindent{\bf Proof.}\ }}
\def\QED{\qed}

\begin{document}
\setlength{\textheight}{7.60truein}    
\runninghead{P. Caldiroli \& R. Musina}{Existence of minimal H-bubbles}
\normalsize\textlineskip
\thispagestyle{empty}
\setcounter{page}{177}
\copyrightheading{Vol.~4, No.~2 (2002) 177--209}
\vspace*{0.88truein}
\fpage{177}
\centerline{\bf EXISTENCE OF MINIMAL H-BUBBLES}
\vspace*{0.37truein}
\centerline{\footnotesize PAOLO CALDIROLI}
\baselineskip=12pt
\centerline{\footnotesize\it Dipartimento di Matematica, 
Universit\`a di Torino,}
\baselineskip=10pt
\centerline{\footnotesize\it via Carlo Alberto, 10 -- 10123 Torino, Italy}
\baselineskip=10pt
\centerline{\footnotesize\it caldiroli@dm.unito.it}
\vspace*{10pt}
\centerline{\footnotesize ROBERTA MUSINA}
\baselineskip=12pt
\centerline{\footnotesize\it Dipartimento di Matematica ed 
Informatica, Universit\`a di Udine,}
\baselineskip=10pt
\centerline{\footnotesize\it via delle Scienze, 206 -- 33100 Udine, Italy}
\baselineskip=10pt
\centerline{\footnotesize\it musina@dimi.uniud.it}
\vspace*{0.225truein}
\publisher{28 June 2000}{4 July 2001}
\vspace*{0.21truein}
\abstracts{Given a function $H\in C^{1}(\R^{3})$ asymptotic to a constant 
at infinity, we investigate the existence of $H$-bubbles, i.e., 
nontrivial, conformal surfaces parametrized by the sphere, with 
mean curvature $H$. Under some global hypotheses we prove the 
existence of $H$-bubbles with minimal energy.}{}{}
\textlineskip                   
\vspace*{14pt}                  
\section{Introduction}      
\noindent
Since 1930, with the renowned papers by Douglas and Rad\'o on 
minimal surfaces, the study of parametric two-dimensional surfaces 
with prescribed mean curvature, satisfying different kinds of 
geometrical or topological side conditions, has constituted a very 
challenging problem and has played a prominent role in the history 
of the Calculus of Variations.

Surfaces with prescribed {\em constant} mean curvature are usually 
known as ``soap films'' or ``soap bubbles''. 
This case has been successfully and deeply investigated by several 
authors, and nowadays a quite wide description of the problem is 
available in the literature (see the survey book by Struwe 
\cite{[Str1]}).

The phenomenon of the formation of an electrified drop is closely 
related to soap film and soap bubbles.
As experimentally observed (see for example \cite{[Gy]}, \cite{[Ga]},
\cite{[Hi]}), an external electric field may affect the shape of the 
drop, and its surface curvature turns out to be nonconstant, in 
general.

However, as regards the mathematical treatment of the case of 
{\em nonconstant} prescribed mean curvature, only few existence 
results of variational type are known.
Apart from few papers on the existence of a ``small'' solution for the 
Plateau problem (we quote, for instance, \cite{[Hi]}, \cite{[Ste1]}
and \cite{[Ste2]}, see also \cite{[CaMu01]}), all the other 
variational-type results hold true in a perturbative setting, namely, 
for curvatures of the form $H(u)=H_{0}+H_{1}(u)$ with
$H_{0}\in\R\setminus\{0\}$ and $H_{1}\in C^{1}(\R^{3})\cap 
L^{\infty}$ having $\|H_{1}\|_{\infty}$ small.
In particular, let us mention the papers \cite{[Str2]}, \cite{[Wa]}, 
\cite{[BeRe]}, \cite{[Ja1]} and \cite{[Ja2]}, which deal with the
Plateau problem, or the corresponding Dirichlet problem.

In this paper we are interested in the existence of 
$\S^{2}$-type parametric surfaces in $\R^{3}$ having prescribed 
mean curvature $H$, briefly, {\em $H$-bubbles}.

More precisely, for $H\in C^{1}(\R^{3})$,
an $H$-bubble is a nonconstant conformal function $\omega\colon 
\R^{2}\to\R^{3}$, smooth as a map on $\S^{2}$, satisfying the 
following problem:
\begin{equation}
    \label{eq:problem}
    \cases{
    \Delta\omega=2H(\omega)\omega_{x}\wedge\omega_{y}&in $\R^{2}$\cr
    \int_{\R^{2}}|\nabla\omega|^{2}<+\infty~.& $ $}
\end{equation}
Here $\omega_{x}=({\partial\omega_{1}\over\partial 
x},{\partial\omega_{2}\over
\partial x},{\partial\omega_{3}\over\partial x})$, 
$\omega_{y}=({\partial
\omega_{1}\over\partial y},{\partial\omega_{2}\over\partial 
y},{\partial
\omega_{3}\over\partial y})$, $\Delta\omega=\omega_{xx}+\omega_{yy}$,
$\nabla\omega=(\omega_{x},\omega_{y})$, and $\wedge$ denotes the 
exterior
product in $\R^{3}$.

In case of nonzero {\em constant} mean curvature $H(u)\equiv H_{0}$,
Brezis and Coron \cite{[BrCo2]} proved that the only nonconstant
solutions to (\ref{eq:problem}) are spheres of radius $|H_{0}|^{-1}$.

In the present paper, we study the existence of $H$-bubbles with 
minimal energy in case $H\colon\R^{3}\to\R$ is a smooth function
satisfying: 
\begin{description}
	\item[$\mathbf{(h_{1})}$]
	$\sup_{u\in\R^{3}}|\nabla H(u+\xi)\cdot u~u|<1$, 
	for some $\xi\in\R^{3}$,
	
	\item[$\mathbf{(h_{\infty})}$]
	$H(u)\to H_{\infty}$ as $|u|\to\infty$, for some
	$H_{\infty}\in\R$.
\end{description}
The assumption $\mathbf{(h_{1})}$ is a global condition on the radial 
component of $\nabla H(\cdot+\xi)$ that, roughly speaking, measures how 
far $H$ differs from a constant.

In addition, we also need that $H$ is nonzero on some sufficiently 
large set. This condition will be made clear in the following.

In order to state our result we need some preliminaries.
Let us point out that problem (\ref{eq:problem}) has a natural 
variational structure, since solutions to (\ref{eq:problem}) are 
formally the critical points of the functional
$$    
\E_{H}(u)={1\over 2}\int_{\R^{2}}|\nabla u|^{2}+
2\int_{\R^{2}}Q(u)\cdot u_{x}\wedge u_{y}~,
$$
where $Q\colon\R^{3}\to\R^{3}$ is any vector field such that
$\div~Q=H$. 

Roughly speaking, the functional $\int_{\R^{2}}Q(u)\cdot u_{x}
\wedge u_{y}$ has the meaning of a volume, for $u$ in a suitable
space of functions.
This is clear when $H(u)\equiv H_{0}$.
Indeed in this case, taking $Q(u)={H_{0}\over 3}u$, one deals with
the standard volume functional $\int_{\R^{2}}u\cdot u_{x}\wedge u_{y}$
which is a determinant homogeneous in $u$ and, for $u$ constant far 
out, measures the algebraic volume enclosed by the surface 
parametrized by $u$.
Moreover, it turns out to be bounded with respect to the Dirichlet 
integral by the Bononcini-Wente isoperimetric inequality.

These facts hold true more generally when $H$ is a bounded nonzero 
function on $\R^{3}$ (see \cite{[Ste1]}). 
In particular, the functional $\int_{\R^{2}}Q(u)\cdot u_{x}
\wedge u_{y}$ is essentially cubic in $u$ and it satisfies a 
generalized isoperimetric inequality.
For this reason, we expect that $\E_{H}$ has a mountain pass 
structure, and this gives an indication for the existence of 
a nontrivial critical point.

The natural space in order to look for $\S^{2}$-type solutions
seems to be the Sobolev space $H^{1}(\S^{2},\R^{3})$, modulo
stereographic projection.
However, working with this space gives some technical difficulties
due to the fact that $H$ may be nonconstant. 
In any case, we can define a mountain pass level for $\E_{H}$ 
restricted to some class of smooth functions. 
In addition, thanks to the assumption $\mathbf{(h_{1})}$, we can 
restrict ourselves to radial paths spanned by functions in
$\Spheres_{\xi}=\{\xi+C^{\infty}_{c}(\R^{2},\R^{3}):
u\not\equiv\xi\}$, where $\xi\in\R^{3}$ is the same as in 
$\mathbf{(h_{1})}$.
Thus we are lead to introduce the value
$$
c_{H}=\inf_{u\in\Spheres_{\xi}}\sup_{s>0}~\E_{H}(su)~.
$$
The assumption $\mathbf{(h_{\infty})}$ guarantees that
$$
0<c_{H}\le\frac{4\pi}{3H_{\infty}^{2}}~.
$$
Note that if $H_{\infty}\ne 0$, the value  
$\frac{4\pi}{3H_{\infty}^{2}}$ equals the mountain pass 
level for the energy functional $\E_{H_{\infty}}$ corresponding to 
the constant mean curvature $H_{\infty}$.
Moreover by the results proved by Brezis and Coron in \cite{[BrCo2]},
this value is the least critical value for $\E_{H_{\infty}}$ in 
$H^{1}(\S^{2},\R^{3})$, and it is attained by the spheres (with 
degree 1) of radius $|H_{\infty}|^{-1}$.
Now, our result can be stated as follows:

\begin{Theorem}
	\label{T:main-result}
	Let $H\in C^{1}(\R^{3})$ satisfy $\mathbf{(h_{1})}$ and
	$\mathbf{(h_{\infty})}$. If 
	$$
	c_{H}<\frac{4\pi}{3H_{\infty}^{2}}
	\leqno(*)
	$$
	holds, then there exists an $H$-bubble $\omega$
	such that $\E_{H}(\omega)=c_{H}$. 
	Moreover, called $\B_{H}$ the set of $H$-bubbles, it holds that 
	$c_{H}=\inf_{\omega\in\B_{H}}\E_{H}(\omega)$.
\end{Theorem}

We point out that, thanks to $\mathbf{(h_{1})}$, the condition
$(*)$ requires that $\E_{H}(\bar u)<0$ for some $\bar u\in
\Spheres_{\xi}$ and then excluded the case $H\equiv 0$.
Clearly, when $\E_{H}(\bar u)<0$ somewhere and $H_{\infty}=0$, then
$(*)$ is automatically satisfied.
Moreover, when $H_{\infty}>0$, the condition $(*)$ turns out to be 
true if $H(u)>H_{\infty}$ for $|u|$ large.
Note that, in general, even if $H(u)=H_{\infty}$ for $|u|\ge R$, 
Theorem \ref{T:main-result} ensures that the $H$-bubble we find is 
different from the $H_{\infty}$-bubble located in the region 
$|u|\ge R$.

We also notice that in general we have no information about the 
position of the $H$-bubble given by Theorem \ref{T:main-result}.
In particular, we can exhibit examples of radial curvatures $H$ for 
which $H$-bubbles with minimal energy exist but cannot be
radial.

The main difficulties in approaching problem (\ref{eq:problem}) with 
variational methods concern the study of the Palais-Smale sequences.
In particular, we emphasize the following problems: boundedness of a 
Palais-Smale sequence with respect to the Dirichlet norm, and in 
$L^{\infty}$; blow up analysis for a (bounded) Palais-Smale sequence.
Concerning the first problem, the assumption $\mathbf{(h_{1})}$ can 
be 
useful in order to guarantee the boundedness with respect to the 
gradient $L^{2}$-norm. However the boundedness in $L^{\infty}$ in 
general cannot be deduced {\em a~priori} and it is not just a 
technical difficulty. In fact, one can exhibit examples of 
Palais-Smale sequences which are bounded with respect to the 
Dirichlet 
norm, but not in $L^{\infty}$, and the lack of boundedness in 
$L^{\infty}$ cannot be eliminated in any way.

Hence, because of these difficulties, we tackle the problem by using 
an approximation method in the spirit of a 
celebrated paper by Sacks and Uhlenbeck \cite{[SaUh]}.
More precisely, we construct a family of approximating solutions 
on which global and local estimates can be proved.
In particular, assuming that $H$ is constant far out, we can obtain 
boundedness both with respect to the Dirichlet norm, and in 
$L^{\infty}$.
Then, a limit procedure, involving a (partial) blow up analysis, is 
carried out, in order to show the existence of an $H$-bubble with 
minimal energy.
In the last step, we remove the assumption that $H$ is constant far 
out, by an approximation argument on the curvature function, and we 
recover the full result stated in Theorem \ref{T:main-result}.

We point out that for a curvature $H\in C^{1}(\R^{3})$ satisfying 
$\mathbf{(h_{1})}$ and such that $H(u)\equiv H_{\infty}\ne 0$ for 
$|u|$ 
large, the set $\B_{H}$ of $H$-bubbles is nonempty {\em a~priori}, 
and the existence of a minimal $H$-bubble can be obtained with a 
direct argument, just minimizing the energy functional $\E_{H}$ 
over $\B_{H}$, without using the above mentioned approximation 
method.
In fact, the hard step lies in removing the 
condition that $H$ is constant far out, just asking to $H$ the 
asymptotic behaviour stated in $\mathbf{(h_{\infty})}$.
To this goal, it is important to know that the energy of the minimal 
$H$-bubble is exactly $c_{H}$, and proving this needs either a sharp 
study of the behaviour of the Palais Smale sequences, or an (almost 
equivalent) approximation argument as, for instance, the 
Sacks-Uhlenbeck type argument that we develop. 
This step requires much more work and constitutes the largest part 
of this paper.

We finally mention a result by Bethuel and Rey \cite{[BeRe]} that 
states the existence of an $H$-bubble passing through an arbitrarily 
prescribed point in $\R^{3}$ in case $H$ is a perturbation of a 
nonzero constant.
This result expresses the 
fact that the bubbles with constant curvature $H_{0}\ne 0$ are 
stable with respect to small $L^{\infty}$ perturbations of $H_{0}$.
Actually, in our opinion, the proof of this result is not completely 
clear and we are not able to recover it with our method.

In fact, we think that the problem of existence of $H$-bubbles for a 
prescribed bounded curvature function $H$ has some similarities with
a semilinear elliptic problem on $\R^{N}$ of the form
\begin{equation}
	\label{eq:elliptic-pb}
	\cases{-\Delta u+u=a(x)u^{p}&on $\R^{N}$\cr
	u>0&on $\R^{N}$\cr u\in H^{1}(\R^{N})}
\end{equation}
where $1<p<\frac{N+2}{N-2}$ and $a$ is a bounded positive function 
on $\R^{N}$.
It is known that the existence of solutions to (\ref{eq:elliptic-pb})
is strongly affected by the behaviour of the coefficient $a(x)$, and 
in some cases 
problem (\ref{eq:elliptic-pb}) has no solution.
In particular, this may happen also when $a(x)$ is a small 
$L^{\infty}$ 
perturbation of a positive constant. 

In our opinion, similar considerations hold also for the problem of 
$H$-bubbles, and the behaviour of $H(u)$ plays a similar role of the
coefficient $a(x)$ in (\ref{eq:elliptic-pb}).
Hence, as well as for problem (\ref{eq:elliptic-pb}), we suspect that 
the existence of $H$-bubbles with minimal energy may depend in a very 
sensitive way on the function $H$.

\section{The variational approach}
\noindent
This Section is structured as follows. In the first part we
introduce some notation in view of setting up a variational 
framework to study problem (\ref{eq:problem}). In particular
we define the $H$-volume functional, the energy functional
associated to problem  (\ref{eq:problem}), and we recall some 
generalized isoperimetric inequality.
In the second part we define a mountain pass level
$c_{H}$ for the energy functional $\E_{H}$ and we 
discuss some properties related to the value $c_{H}$
strongly depending on the assumption $\mathbf{(h_{1})}$. 

\subsection{Notation and isoperimetric inequality}
\noindent
First, let us introduce the space 
$$
X=\{v\circ\phi:v\in H^{1}(\S^{2},\R^{3})\}
$$
where $\phi\colon\R^{2}\to\S^{2}$ is the (inverse of the) standard 
stereographic projection and it is given by
\begin{equation}
	\phi(z)=(\mu x,\mu y,1-\mu)~,
	\ \ \mu=\mu(z)=\frac{2}{1+|z|^{2}}~,
	\label{eq:stereographic-projection}
\end{equation}
being $z=(x,y)$ and $|z|^{2}=x^{2}+y^{2}$.
Notice that $u\in X$ if and only if $u,\hat u\in H^{1}_{loc}
(\R^{2},\R^{3})$ and $\int_{\R^{2}}|\nabla u|^{2}<+\infty$, where
$\hat u(z)=u\big(\frac{z}{|z|^{2}}\big)$.
Let us also set $H^{1}_{0}=H^{1}_{0}(\disc,\R^{3})$, where $\disc$ 
is the open unit disc in $\R^{2}$. Clearly, $H^{1}_{0}\subset X$.
For every $u\in X$ we denote the Dirichlet integral by
$$
\D(u)=\frac{1}{2}\int_{\R^{2}}|\nabla u|^{2}~.
$$

Now, given $H\in C^{1}(\R^{3})$, we construct the $H$-volume 
functional as follows.
Set 
$$
m_{H}(u)=\int_{0}^{1}H(su)s^{2}~ds~.
$$
Thus, for every $u\in\R^{3}$ one has
\begin{equation}
	\label{eq:div}
	\div(m_{H}(u)u)=H(u)~.
\end{equation}
Then, let $\V_{H}\colon X\cap L^{\infty}\to\R$ be defined by
$$
\V_{H}(u)=\int_{\R^{2}}m_{H}(u)u\cdot u_{x}\wedge u_{y}~.
$$
In case $H(u)\equiv 1$, one has $m_{H}(u)\equiv\frac{1}{3}$, and the 
functional $\V_{H}$ reduces to the classical volume functional which 
satisfies the standard isoperimetric inequality.
In fact the following generalization holds, as proved by Steffen in
\cite{[Ste1]}.

\begin{Lemma}
	\label{L:isoperimetric-inequality}
	If $H\in C^{1}(\R^{3})$ is bounded on $\R^{3}$ then 
	there exists $S_{H}>0$ such that
	\begin{equation}
		\label{eq:isoperimetric-inequality}
		S_{H}|\V_{H}(u)|^{2/3}\le\D(u)\ \ 
		{\it for\ every\ }u\in X\cap L^{\infty}~.
	\end{equation}
\end{Lemma}

\begin{Remark}
	\label{R:isoperimetric-inequality}
	{\rm In fact Steffen in \cite{[Ste1]} proves that the 
	functional $\V_{H}$ admits a continuous extension on 
	$H^{1}_{0}$ and (\ref{eq:isoperimetric-inequality}) 
	holds true also for every $u\in H^{1}_{0}$.}
\end{Remark}

Finally we introduce the energy functional 
$\E_{H}\colon X\cap L^{\infty}\to\R$, defined for every 
$u\in X\cap L^{\infty}$ by
$$
\E_{H}(u)=\D(u)+2\V_{H}(u)~.
$$
In the following result we state some properties of the functional 
$\E_{H}$. 

\begin{Lemma}
	\label{L:functional}
	Let $H\in C^{1}(\R^{3})$.
	Then:
	\begin{description}
		\item[$(i)$]
		~~for every $u\in X\cap L^{\infty}$ one has 
		$\E_{H}(su)=s^{2}\D(u)+o(s^{2})$ as $s\to 0$,
		
		\item[$(ii)$]
		~~for every $u\in X\cap L^{\infty}$ and for $h\in 
		C^{\infty}_{c}(\R^{2},\R^{3})$ the directional derivative
		of $\E_{H}$ at $u$ along $h$ exists, and it is given by
		$$
		{d\E_{H}(u)h}=\int_{\R^{2}}\nabla u\cdot\nabla h
		+2\int_{\R^{2}}H(u) h\cdot u_{x}\wedge u_{y}~,
		$$
		
		\item[$(iii)$]
		~~for every bounded solution $\omega$ to 
		{\rm (\ref{eq:problem})} one has 
		\begin{equation}
			\D(\omega)+\int_{\R^{2}}H(\omega)
			\omega\cdot\omega_{x}\wedge\omega_{y}=0~.
			\label{eq:omega-condition}
		\end{equation}
	\end{description}
\end{Lemma}

\begin{Remark}
	\label{R:H-bubble-regularity}
	{\rm If $\omega\in X\cap L^{\infty}$ is a weak solution to 
	(\ref{eq:problem}), i.e., $d\E_{H}(\omega)h=0$ for every 
	$h\in C^{\infty}_{c}(\R^{2},\R^{3})$, then, since $H\in 
	C^{1}(\R^{3})$, by a Heinz regularity result \cite{[He]}, 
	$\omega\in C^{3}(\R^{2},\R^{3})$, it is conformal, and smooth 
	as a map on $\S^{2}$. 
	In particular there exists 
	$\lim_{|z|\to\infty}\omega(z)=\omega_{\infty}\in\R^{3}$.}   
\end{Remark}

\Proof
Part $(i)$ is a consequence of Lemma \ref{L:isoperimetric-inequality}.
Part $(ii)$ follows by the results in \cite{[HiKa]}, 
using (\ref{eq:div}). Finally, (\ref{eq:omega-condition}) can be 
proved multiplying the system $\Delta\omega=2H(\omega)\omega_{x}
\wedge\omega_{y}$ by $\omega$, integrating on $\disc_{R}$, and 
passing to the limit as $R\to+\infty$.
\QED

To conclude this Subsection, we point out a consequence of 
assumption $\mathbf{(h_{\infty})}$. Actually, the following 
result holds true under a much weaker condition.

\begin{Lemma}
	\label{L:negative-value}
	Let $H\in C^{1}(\R^{3},\R)$ satisfy
	\begin{equation}
		|H(su)|\ge H_{0}>0\ \ {\it for}\ s\ge s_{0}
		\ {\it and}\ u\in\Sigma~,
		\label{eq:negative-value}
	\end{equation}
	being $\Sigma$ a nonempty open set in $\S^{2}$.
	Then there exists $\bar u\in H^{1}_{0}\cap L^{\infty}$
	such that $\E_{H}(s\bar u)\to-\infty$ as $s\to+\infty$.
\end{Lemma}

\Proof
Thanks to the rotational invariance of the problem we may assume 
that $\Sigma$ is an open neighborhood of the point $-e_{3}=(0,0,-1)$.
Furthermore, let us suppose that $H(su)\ge H_{0}>0$ for $s>s_{0}$
and $u\in\Sigma$.
For $\delta\in(0,1)$ let us define 
$$
u^{\delta}(z)=\cases{
\phi(z)& as $|z|<\delta$\cr
\frac{1-|z|}{1-\delta}~\phi\left(\frac{\delta}{|z|}z\right)
& as $\delta\le|z|\le 1~,$}
$$
where $\phi\colon\R^{2}\to\S^{2}$ is the function introduced in
(\ref{eq:stereographic-projection}).
It holds that $u^{\delta}\in H^{1}_{0}\cap L^{\infty}$, and 
$u^{\delta}$ parametrizes the boundary of the sector of cone 
defined by
$$
A_{\delta}=\{\xi\in\R^{3}:-|\xi|\cos\theta_{\delta}>
\xi\cdot e_{3}~,\ |\xi|<1\}
$$
where $\theta_{\delta}=\arccos\frac{1-\delta^{2}}{1+\delta^{2}}$.
In addition one has that
$$
u^{\delta}(z)\cdot u^{\delta}_{x}(z)\wedge u^{\delta}_{y}(z)=
\cases{-\mu(z)^{2}& as $|z|<\delta$\cr
0& as $|z|>\delta$}
$$
and for every $s>0$, by the divergence theorem,
$$
\V_{H}(su^{\delta})=-\int_{sA_{\delta}}H(\xi)~d\xi~.
$$
Since $\phi(0)=-e_{3}$ and $\phi$ is continuous, we can 
find $\delta_{0}\in(0,1)$ such that $\phi(z)\in\Sigma$ as 
$|z|<\delta_{0}$.
Set $\bar u=u^{\delta_{0}}$ and $A=A_{\delta_{0}}$.
Therefore, by the hypothesis, for $s>s_{0}$ one has
\begin{eqnarray*}
	\V_{H}(s\bar u)&=&-\int_{s_{0}A}H(\xi)~d\xi
	-\int_{sA\setminus s_{0}A}H(\xi)~d\xi\\
	&\le&\V_{H}(s_{0}\bar u)-\int_{sA\setminus s_{0}A}H_{0}~d\xi\\
	&=&\V_{H}(s_{0}\bar u)-\V_{H_{0}}
	(s_{0}\bar u)-s^{3}|\V_{H_{0}}(\bar u)|~.
\end{eqnarray*}
Then
$$
\E_{H}(s\bar u)\le s^{2}\D(\bar u)+2(\V_{H}(s_{0}\bar u)-
\V_{H_{0}}(s_{0}\bar u))-2s^{3}|\V_{H_{0}}(\bar u)|~.
$$
Passing to the limit as $s\to+\infty$ we obtain the thesis.
Finally, we observe that in case $H(su)\le H_{0}<0$ for
$s>s_{0}$ and $u\in\Sigma$, one can repeat the same argument 
taking $v(x,y)=u(y,x)$.
\QED

\subsection{The mountain pass level}
\noindent
Assume that $H\in 
C^{1}(\R^{3})\cap L^{\infty}$ is such that  
there exists $\bar u\in C^{\infty}_{c}(\R^{2},\R^{3})$ with
$\E_{H}(\bar u)<0$.
In particular, this excludes the case $H\equiv 0$. 
Then, let
\begin{equation}
	\label{eq:radial-mp-level}
	c_{H}=\inf_{u\in C^{\infty}_{c}(\R^{2},\R^{3})\atop u\ne 0} 
	\sup_{s>0}~\E_{H}(su)~.
\end{equation} 
Note that $c_{H}$ is well defined and, thanks to Lemma 
\ref{L:isoperimetric-inequality}, it is positive and finite.
In particular, by
(\ref{eq:isoperimetric-inequality}), one can estimate 
$$
c_{H}\ge\left(\frac{S_{H}}{3}\right)^{3},
$$
where $S_{H}$ is the isoperimetric constant 
associated to $H$.

\begin{Remark}
	\label{R:constant-curvature}
	{\rm When $H(u)\equiv H_{0}\in\R\setminus\{0\}$, the volume 
	functional is purely cubic and one can easily prove that 
	$$
	c_{H_{0}}=
	\left(\frac{S_{H_{0}}}{3}\right)^{3}=
	\frac{4\pi}{3H_{0}^{2}}=\E_{H_{0}}(\omega^{0})
	=\sup_{s>0}\E_{H_{0}}(s\omega^{0})
	$$ 
	where $\omega^{0}=\frac{1}{H_{0}}\phi$ and $\phi$ is defined 
	in (\ref{eq:stereographic-projection}).
	Notice that $\omega^{0}$ is a conformal parametrization of the  
	sphere of radius $|H_{0}|^{-1}$ centered at the origin, 
	it satisfies $\Delta\omega^{0}=2H_{0}\omega^{0}_{x}
	\wedge\omega^{0}_{y}$ on $\R^{2}$, 
	$\D(\omega^{0})=\frac{4\pi}{H_{0}^{2}}$, and 
	$\V_{H_{0}}(\omega^{0})=-\frac{4\pi}{3H_{0}^{2}}$.}
\end{Remark}

The results that follow better explain the role of the condition 
$\mathbf{(h_{1})}$ with respect to the definition of $c_{H}$. 
To this extent, we point out that, since problem (\ref{eq:problem}) 
is invariant under translations, in the assumption 
$\mathbf{(h_{1})}$ we may suppose that $\xi=0$.
Hence, setting:
\begin{equation}
	\label{eq:MH-definition}
	M_{H}=\sup_{u\in\R^{3}}|\nabla H(u)\cdot u~u|
\end{equation}
the hypothesis $\mathbf{(h_{1})}$ reads: $M_{H}<1$. 
It is convenient to introduce also the value
\begin{equation}
	\label{eq:barMH-definition}
	\bar M_{H}=2\sup_{u\in\R^{3}}|(H(u)-3m_{H}(u))u|~.
\end{equation}
In fact, several estimates in the sequel need a bound just on
$\bar M_{H}$.

\begin{Remark}
	\label{R:MH}
	{\rm $(i)$ By (\ref{eq:div}) and by the definition of $m_{H}$, 
	it turns out that $\bar M_{H}\le M_{H}$, but the strict inequality 
	may also occur. Indeed one can construct functions $H\in C^{1}
	(\R^{3})$ such that $M_{H}=+\infty$, while $\bar M_{H}<+\infty$.	
	
	\noindent
	$(ii)$ If $H\in C^{1}(\R^{3})$ satisfies $\bar M_{H}<+\infty$
	then it turns out that $H\in L^{\infty}(\R^{3})$. Furthermore, 
	for every $u\in\S^{2}$ there exists $\lim_{s\to+\infty}H(su)=
	\hat H(u)\in\R$	and $\hat H\in C^{0}(\S^{2})$. Thus, if 
	$\bar M_{H}<+\infty$, then 
	the condition (\ref{eq:negative-value}) used in Lemma 
	\ref{L:negative-value} is verified whenever 
	$\limsup_{s\to+\infty}|H(su)|>0$ for some $u\in\S^{2}$.
	}
\end{Remark}

First, we give a positive lower bound on the
energy of any $H$-bubble.

\begin{Proposition}
	\label{P:minimal-energy}
	Let $H\in C^{1}(\R^{3})$ satisfy $\mathbf{(h_{1})}$.
	If $\omega$ is an $H$-bubble, then $\E_{H}(\omega)\ge
	{c}_{H}$.
\end{Proposition}

\noindent
The proof of Proposition \ref{P:minimal-energy} is based
on the following Lemma.

\begin{Lemma}
	\label{L:radial-mp-level}
	Let $H\in C^{1}(\R^{3})$ satisfy $\bar M_{H}<1$ and let
	$u\in H^{1}_{0}\cap L^{\infty}\setminus\{0\}$.
	\begin{description}
		\item[$(i)$] 
		If $\sup_{s>0}\E_{H}(su)<+\infty$ then
		$\E_{H}(su)\le as^{2}-bs^{3}$ for every $s>0$, with $a,b>0$
		depending on $u$, 
		
		\item[$(ii)$] 
		if $\E_{H}(s_{0}u)<0$ for some $s_{0}>0$ then 
		$\sup_{s>0}\E_{H}(su)=\max_{s\in[0,s_{0}]}\E_{H}(su)$,
		
		\item[$(iii)$] 
		if $\sup_{s>0}\E_{H}(su)=\E_{H}(\bar su)$,
		then $\V_{H}(\bar su)<0$.

	\end{description}
\end{Lemma}

\Proof
Fix $u\in H^{1}_{0}\cap L^{\infty}\setminus\{0\}$ and set $f(s)=
\E_{H}(su)$ for every $s\ge 0$.
Notice that $f$ is differentiable and 
$$
f'(s)=
s\int_{\disc}|\nabla u|^{2}+2s^{2}
\int_{\disc}H(su)u\cdot u_{x}\wedge u_{y}~.
$$
Using (\ref{eq:barMH-definition}), one has that
\begin{equation}
	\label{eq:f1-inequality}
	f'(s)\le-(1-\bar M_{H})\D(u)s+\frac{3}{s}f(s)~.
\end{equation}
If $\sup_{s>0}f(s)<+\infty$, since $\bar M_{H}<1$,
(\ref{eq:f1-inequality}) implies
that $\lim_{s\to+\infty}f'(s)=-\infty$ and then there exists 
$s_{0}>0$ such that $f(s)<0$ for $s\ge s_{0}$.
Setting $\bar a=(1-\bar M_{H})\D(u)$ and integrating 
(\ref{eq:f1-inequality}) over $[s_{0},s]$ one obtains
\begin{equation}
	\label{eq:f1-integrated}
	f(s)\le\left(\frac{f(s_{0})}{s_{0}^{3}}-\frac{\bar a}
	{s_{0}}\right)s^{3}+\bar as^{2}
\end{equation}
for every $s\ge s_{0}$.
Keeping into account that $f(s)=s^{2}\D(u)+o(s^{2})$ 
as $s\to 0^{+}$, one can find $a\ge\bar a$ such that
$$
f(s)\le\left(\frac{f(s_{0})}{s_{0}^{3}}-\frac{\bar a}
{s_{0}}\right)s^{3}+as^{2}
$$
for every $s\ge 0$, namely $(i)$.
Now, let us prove $(ii)$. If $f(s_{0})<0$, by 
(\ref{eq:f1-integrated}), one infers that $\sup_{s>0}f(s)<+\infty$.
Moreover (\ref{eq:f1-inequality}) implies in particular that 
$f'(s)<0$ whenever $f(s)\le 0$. Hence also $(ii)$ holds true.
Finally, if $\sup_{s>0}\E_{H}(su)=\E_{H}(\bar su)$,
then $f'(\bar s)=0$, 
and consequently, by (\ref{eq:barMH-definition}),
$$
3\V_{H}(\bar su)=3\V_{H}(\bar su)-\bar sf'(\bar s)
\le-\bar s^{2}\left(1-\frac{\bar M_{H}}{2}\right)
\D(u)<0~, 
$$
that is $(iii)$.
\QED

\noindent
{\bf Proof of Proposition \ref{P:minimal-energy}.}
By Remark \ref{R:H-bubble-regularity} an $H$-bubble $\omega$ is 
smooth and bounded.
Moreover the mapping $f(s)=\E_{H}(s\omega)$ is well defined, and
twice differentiable on $(0,+\infty)$, with
$$
f''(s)=\int_{\disc}|\nabla\omega|^{2}+4s\int_{\disc}H(s\omega)
\omega\cdot\omega_{x}\wedge\omega_{y}+2s^{2}\int_{\disc}
\nabla H(s\omega)\cdot\omega~\omega\cdot\omega_{x}\wedge
\omega_{y}~.
$$
Since $M_{H}<1$, one obtains that
\begin{equation}
	\label{eq:f2-inequality}
	f''(s)\le-2(1-M_{H})\D(u)+\frac{2}{s}f'(s)~.
\end{equation}
In particular, by (\ref{eq:f2-inequality}), if $f'(\bar s)=0$ for 
some $\bar s>0$ then $f''(\bar s)<0$. This shows that there exists 
at most one value $\bar s>0$ where $f'(\bar s)=0$.
In fact, one knows that $f'(1)=0$ because of 
(\ref{eq:omega-condition}). Hence $\sup_{s>0}f(s)=f(1)$ and, 
arguing as in the proof of Lemma \ref{L:radial-mp-level}, 
$f(s)\to-\infty$ as $s\to+\infty$.
Now, for every $\delta\in(0,1)$ let $u^{\delta}\colon\disc
\to\R^{3}$ be defined as follows:
$$
u^{\delta}(z)=
\cases{
0 & as $|z|\ge{\delta}$\cr
\big(\frac{\log|z|}{\log\delta}-1\big)
\omega_{\infty} & as $\delta^{2}\le |z|<\delta$\cr
\big(\frac{\log|z|}{2\log\delta}-1\big)
(\omega^{\delta}(z)-\omega_{\infty})
+\omega_{\infty}& as $\delta^{4}\le |z|<\delta^{2}$\cr
\omega^{\delta}(z) & as $|z|<\delta^{4}$}
$$
where $\omega_{\infty}=\lim_{|z|\to\infty}\omega(z)$,
and $\omega^{\delta}(z)=\omega({z\over\delta^{5}})$.
Note that $u^{\delta}\in H^{1}_{0}\cap L^{\infty}$ and
$\|u^{\delta}\|_{\infty}\le\|\omega\|_{\infty}$.
Let us set $f_{\delta}(s)=\E_{H}(su^{\delta})$.
We claim that for every $s'>0$
\begin{equation}
	\sup_{s\in[0,s']}|f_{\delta}(s)-f(s)|\to 0\quad
	{\rm as}\ \delta\to 0~.
	\label{eq:claim}
\end{equation}
Assuming for a moment that (\ref{eq:claim}) holds, let us complete 
the proof.
Let $s_{0}>1$ be such that $f(s_{0})<0$. 
By (\ref{eq:claim}), for $\delta>0$ small enough, 
$f_{\delta}(s_{0})<0$ and then, by Lemma \ref{L:radial-mp-level},
$\sup_{s>0}f_{\delta}(s)$ is attained in $(0,s_{0})$.
Hence, using again (\ref{eq:claim}), we have
$$
{c}_{H}\le\sup_{s>0}f_{\delta}(s)=
\max_{s\in[0,s_{0}]}f_{\delta}(s)\le
\max_{s\in[0,s_{0}]}f(s)+o(1)=f(1)+o(1)~.
$$
Therefore the thesis follows.
Finally, let us prove the claim (\ref{eq:claim}).
For every $s\ge 0$ we can write
\begin{eqnarray*}
	f_{\delta}(s)-f(s)&=&
	\!s^{2}\left(\int_{|z|>\delta^{4}}|\nabla u^{\delta}|^{2}-
	\int_{|z|>\delta^{-1}}|\nabla\omega|^{2}\right)\\
	&+&\!2s^{3}\left(\int_{|z|>\delta^{4}}\!m_{H}(su^{\delta})
	u^{\delta}\cdot u^{\delta}_{x}\wedge u^{\delta}_{y}-
	\int_{|z|>\delta^{-1}}\!m_{H}(s\omega)\omega\cdot
	\omega_{x}\wedge\omega_{y}\right)\!.
\end{eqnarray*}
We observe that
\begin{eqnarray*}
	2\left|\int_{|z|>\delta^{4}}m_{H}(su^{\delta})u^{\delta}\cdot
	u^{\delta}_{x}\wedge u^{\delta}_{y}\right|&\le&
	\|m_{H}\|_{\infty}\|\omega\|_{\infty}
	\int_{|z|>\delta^{4}}|\nabla u^{\delta}|^{2}\\
	2\left|\int_{|z|>\delta^{-1}}m_{H}(s\omega)\omega\cdot
	\omega_{x}\wedge\omega_{y}\right|&\le&
	\|m_{H}\|_{\infty}\|\omega\|_{\infty}
	\int_{|z|>\delta^{-1}}|\nabla\omega|^{2}.
\end{eqnarray*}
Moreover, one can check that
$$
\int_{|z|>\delta^{4}}|\nabla u^{\delta}|^{2}\to 0\quad
{\rm as}\ \delta\to 0~,
$$
and, since $\omega\in X$, also
$$
\int_{|z|>\delta^{-1}}|\nabla\omega|^{2}\to 0\quad
{\rm as}\ \delta\to 0~.
$$
Therefore (\ref{eq:claim}) immediately follows and this concludes 
the proof.
\QED

Notice that the full condition $M_{H}<1$ enters just in the previous 
step.
Now we are going to prove two technical Lemmata that will be used
in the sequel.

\begin{Lemma}
	\label{L:c-lambda-H}
	Let $H\in C^{1}(\R^{3})$ satisfy $\bar M_{H}<1$.
	Then $c_{H}\le c_{\lambda H}$ 
	for every $\lambda\in(0,1]$.
\end{Lemma}

\Proof
Firstly, notice that for $\lambda\in(0,1]$, the isoperimetric 
inequality 
(\ref{eq:isoperimetric-inequality}) holds true also for 
$\lambda H$ (with $S_{\lambda H}=\lambda^{-\frac{2}{3}}S_{H}$),
and then the value $c_{\lambda H}$ is well defined and positive.
Suppose that it is finite and, given $\epsilon>0$, let $u\in 
C^{\infty}_{c}(\R^{2},\R^{3})\setminus\{0\}$ be such that 
$\sup_{s>0}\E_{\lambda H}(su)<c_{\lambda H}+\epsilon$.
Since $\bar M_{\lambda H}=\lambda\bar M_{H}<1$, by Lemma 
\ref{L:radial-mp-level}, $\lim_{s\to+\infty}\E_{\lambda H}(su)
=-\infty$. In particular, $\V_{\lambda H}(su)<0$ for $s$ large.
Hence $\E_{H}(su)\le \E_{\lambda H}(su)<0$ for $s$ large.
Using again Lemma \ref{L:radial-mp-level} there exists  
$\bar s>0$ such that $\sup_{s>0}\E_{H}(su)=\E_{H}(\bar su)$. 
Furthermore $\V_{H}(\bar su)<0$. 
Therefore
$$
c_{H}\le\E_{H}(\bar su)=
\E_{\lambda H}(\bar su)+2(1-\lambda)\V_{H}(\bar su)\le
\E_{\lambda H}(\bar su)\le\sup_{s>0}\E_{\lambda H}(su)
\le c_{\lambda H}+\epsilon~.
$$
Then the thesis 
follows because of the arbitrariness of $\epsilon>0$.
\QED

The next result states the upper semicontinuity of 
$c_{H}$ with respect to $H$.

\begin{Lemma}
	\label{L:semicontinuity-c-H}
	Let $H\in C^{1}(\R^{3})$ satisfy $\bar M_{H}<1$. Let 
	$(H_{n})\subset C^{1}(\R^{3})$ be a sequence of functions 
	satisfying $\bar M_{H_{n}}<1$, and such that $H_{n}\to H$ 
	uniformly on compact sets of $\R^{3}$.
	Then $\limsup_{n\to+\infty}c_{H_{n}}
	\le{c}_{H}$.	
\end{Lemma}

\Proof
Suppose that ${c}_{H}$ is finite and, given $\epsilon>0$ 
take $u\in C^{\infty}_{c}(\R^{2},\R^{3})\setminus\{0\}$ such that
$\sup_{s>0}\E_{H}(su)<{c}_{H}+\epsilon$.
One can check that $\lim_{n\to+\infty}\E_{H_{n}}(su)=\E_{H}(su)$ 
for every $s\ge 0$. By Lemma \ref{L:radial-mp-level}, 
$\E_{H}(s_{0}u)<0$ for some $s_{0}>0$, and then also 
$\E_{H_{n}}(s_{0}u)<0$ for $n\in\N$ large enough. 
Therefore, since $H_{n}$ satisfies $\bar M_{H_{n}}<1$, using again 
Lemma \ref{L:radial-mp-level}, $\sup_{s>0}\E_{H_{n}}(su)=
\E_{H_{n}}(\bar s_{n}u)$ for some $\bar s_{n}\in[0,s_{0}]$. 
Then, for a subsequence, $\bar s_{n}\to\bar s$ and, since 
$H_{n}\to H$ uniformly on compact sets, $\E_{H_{n}}(\bar s_{n}u)\to
\E_{H}(\bar s u)$. Consequently one has
$$
{c}_{H_{n}}\le\E_{H_{n}}(\bar s_{n}u)=
\E_{H}(\bar s u)+o(1)\le\sup_{s>0}\E_{H}(su)
+o(1)\le{c}_{H}+\epsilon+o(1).
$$
Passing to the limit as $n\to+\infty$ and taking into account of the 
arbitrariness of $\epsilon>0$, the thesis is proved.
\QED

Lastly, we give an estimate for
${c}_{H}$ from above.
Here, just the assumption $\mathbf{(h_{\infty})}$, and in fact
a more general condition, is enough.

\begin{Lemma}
	\label{L:large-inequality}
	Let $H\in C^{1}(\R^{3})$ satisfy 
	{\rm (\ref{eq:negative-value})} for some nonempty open set 
	$\Sigma\subset\S^{2}$.
	Then ${c}_{H}\le\frac{4\pi}{3H_{0}^{2}}$.
\end{Lemma}

\Proof
As in the proof of Lemma \ref{L:negative-value}, we may assume 
that $\Sigma$ is an open neighborhood of the point $-e_{3}=(0,0,-1)$
and that $H(su)\ge H_{0}>0$ for $s>s_{0}$ and $u\in\Sigma$.
Let us consider the function $\omega^{0}\colon\R^{2}\to\R^{3}$ 
defined 
as in Remark \ref{R:constant-curvature}. 
For every $r>0$ set $\omega^{r}=\omega^{0}-re_{3}$.
Notice that $\omega^{r}$ is a conformal parametrization of a sphere
of radius $r_{0}=\frac{1}{H_{0}}$ and center $-re_{3}$.
Hence, using the divergence theorem, one has that
\begin{equation}
	\V_{H}(s\omega^{r})=-\int_{B_{sr_{0}}(-se_{3})}H(\xi)~d\xi=
	-s^{3}\int_{B_{r_{0}}(0)}H(s\xi-se_{3})~d\xi~.
	\label{eq:large-inequality-1}
\end{equation}
Setting $s_{r}=\frac{s_{0}}{r-r_{0}}$ and using 
(\ref{eq:large-inequality-1}), one obtains that for $s\in[0,s_{r}]$ 
$$
\E_{H}(s\omega^{r})\le 4\pi(r_{0}s_{r})^{2}
+\frac{8\pi}{3}\|H\|_{\infty}(r_{0}s_{r})^{3}=O(s_{r}^{2})~,
$$
while, for $s\ge s_{r}$, by
the hypothesis (\ref{eq:negative-value}), one has
$$
\E_{H}(s\omega^{r})\le 4\pi(r_{0}s)^{2}
-\frac{8\pi}{3}H_{0}(r_{0}s)^{3}\le\frac{4\pi}{3H_{0}^{2}}~.
$$
Then
\begin{equation}
	\sup_{s>0}\E_{H}(s\omega^{r})\le\max\left\{\frac{4\pi}{3H_{0}^{2}},
	O(s_{r}^{2})\right\}=\frac{4\pi}{3H_{0}^{2}}
	\label{eq:large-inequality-2}
\end{equation}
for $r>0$ large enough.
Now, as in the proof of Proposition \ref{P:minimal-energy}, one can 
construct $u^{r,\delta}\in H^{1}_{0}\cap L^{\infty}$ such that 
$\sup_{s>0}\E_{H}(su^{r,\delta})\le\sup_{s>0}\E_{H}(s\omega^{r})
+o(1)$, with $o(1)\to 0$ as $\delta\to 0$.
Hence, by (\ref{eq:large-inequality-2}), one obtains $c_{H}\le
\frac{4\pi}{3H_{0}^{2}}+o(1)$, that is, the thesis.
\QED

{From} the previous proof, one immediately infers the next estimate.

\begin{Corollary}
	\label{C:strong-inequality}
	Let $H\in C^{1}(\R^{3})$ satisfy 
	$\mathbf{(h_{\infty})}$.
	Then ${c}_{H}\le\frac{4\pi}{3H_{\infty}^{2}}$.
	If, in addition, $H(u)>H_{\infty}>0$ for $|u|$ large, then
	${c}_{H}<\frac{4\pi}{3H_{\infty}^{2}}$.
\end{Corollary}

\section{Approximating problems}
\noindent
Aim of this Section is to introduce a family of perturbed energy
functionals having a mountain pass critical point at a level which
approximate the value $c_{H}$ introduced in the previous Section.

The advantage in following this procedure (already used in a 
different framework by Sacks and Uhlenbeck \cite{[SaUh]}) is due to 
the possibility to obtain some uniform global and local estimates on 
the critical points of the perturbed problems.

Thus, for every $\alpha>1$ ($\alpha$ will be taken close to 1)
we consider the Sobolev space 
$H^{1,2\alpha}_{0}=H^{1,2\alpha}_{0}(\disc,\R^{2})$ and the functional
$\E^{\alpha}_{H}\colon H^{1,2\alpha}_{0}\to\R$ defined by
$$
\E^{\alpha}_{H}(u)=
{1\over 2\alpha}\int_{\disc}\left((1+|\nabla u|^{2})^{\alpha}-1\right)
+2\V_{H}(u)~.
$$
It is convenient to denote
$$
\D^{\alpha}(u)={1\over 2\alpha}\int_{\disc}
\left((1+|\nabla u|^{2})^{\alpha}-1\right)~.
$$
Since $H^{1,2\alpha}\hookrightarrow L^{\infty}\cap H^{1}$, the 
functional $\E^{\alpha}_{H}$ turns out to be well defined and regular 
on $H^{1,2\alpha}_{0}$, when $H$ is any bounded, smooth function. 
More precisely, $\E^{\alpha}_{H}$ is of class $C^{1}$ on 
$H^{1,2\alpha}_{0}$ and
$$
d\E^{\alpha}_{H}(u)h=\int_{\disc}(1+|\nabla u|^{2})^{\alpha-1}
\nabla u\cdot\nabla h+2\int_{\disc}H(u)h\cdot u_{x}\wedge u_{y}
$$
for every $u,h\in H^{1,2\alpha}_{0}$ (see \cite{[HiKa]}).

Our first goal is to prove that for every $\alpha>1$ sufficiently 
close to 1 the functional $\E^{\alpha}_{H}$ has a mountain pass 
geometry and a corresponding mountain pass critical point, as stated
in the following result.

\begin{Lemma}
	\label{L:approximating-solutions}
	Let $H\in C^{1}(\R^{3})\cap L^{\infty}$ be such that 
	there exists $\bar u\in
	C^{\infty}_{c}(\disc,\R^{3})$ with $\E_{H}(\bar u)<0$. 
	Then there exists $\bar\alpha>1$ such that for every 
	$\alpha\in(1,\bar\alpha)$ the class $\Gamma^{\alpha}=
	\{\gamma\in C([0,1],H^{1,2\alpha}_{0}):\gamma(0)=0,~
	\E^{\alpha}_{H}(\gamma(1))<0\}$ is nonempty and the value
	$$
	\overline{c}^{\alpha}_{H}=
	\inf_{\gamma\in\Gamma^{\alpha}}\max_{s\in[0,1]}
	\E^{\alpha}_{H}(\gamma(s))
	$$
	is positive.
	
	\noindent
	If in addition $\bar M_{H}<+\infty$ then for every 
	$\alpha\in(1,\bar\alpha)$ there exists $u^{\alpha}\in 
	H^{1,2\alpha}_{0}$ such 
	that $\E^{\alpha}_{H}(u^{\alpha})=\overline{c}^{\alpha}_{H}$ 
	and $d\E^{\alpha}_{H}(u^{\alpha})=0$.
\end{Lemma}

The second step consists in obtaining some uniform 
estimates on the mountain pass critical points $u^{\alpha}$ 
of the perturbed functionals $\E_{H}^{\alpha}$.

\begin{Proposition}
	\label{P:global-estimates}
	Let $H\in C^{1}(\R^{3})$ be such that $\bar M_{H}<1$ and,
	for every $\alpha\in(1,\bar\alpha)$, let $u^{\alpha}\in 
	H^{1,2\alpha}$ be the critical point of $\E_{H}^{\alpha}$ 
	at level $\overline{c}_{H}^{\alpha}$ given by Lemma 
	\ref{L:approximating-solutions}.
	Then
	\begin{eqnarray*}
		& &\limsup_{\alpha\to 1}\E_{H}^{\alpha}(u^{\alpha})
		\le{c}_{H}~,\\
		& &\sup_{\alpha\in(1,\bar\alpha)}
		\|\nabla u^{\alpha}\|_{2}<+\infty~,\\
		& &\inf_{\alpha\in(1,\bar\alpha)}
		\|\nabla u^{\alpha}\|_{2}>0~,
	\end{eqnarray*}
	where ${c}_{H}$ is defined by {\rm (\ref{eq:radial-mp-level})}.
	If, in addition, $H(u)=H_{0}$ for $|u|\ge R_{0}$, for some 
	$R_{0}>0$, then
	$$
	\sup_{\alpha\in(1,\bar\alpha)}
	\|u^{\alpha}\|_{\infty}<+\infty~.
	$$
\end{Proposition}

The proofs of Lemma \ref{L:approximating-solutions} and 
Proposition
\ref{P:global-estimates} will be carried out in Subsections
3.1 and 3.2, respectively.

The last result of this Section states the behaviour of the
family of the mountain pass critical points $u^{\alpha}$ 
in the limit as $\alpha\to 1$.
This result describes a blow up phenomenon, and it will be proved
in the Appendix, in a more general situation.

\begin{Proposition}
	\label{P:blow-up-0}
	Let $H\in C^{1}(\R^{3})$ be such that $\bar M_{H}<1$ and
	$H(u)=H_{0}$ for $|u|\ge R_{0}$, for some $R_{0}>0$.
	For every $\alpha\in(1,\bar\alpha)$, let $u^{\alpha}\in 
	H^{1,2\alpha}$ be the critical point of $\E_{H}^{\alpha}$ 
	at level $\overline{c}_{H}^{\alpha}$ given by Lemma 
	\ref{L:approximating-solutions}.
	Then, there exist sequences 
	$(\epsilon_{\alpha})\subset(0,+\infty)$, 
	$(z_{\alpha})\subset\overline{\disc}$,
	a number $\lambda\in(0,1]$, and a function 
	$\omega\in X\cap L^{\infty}$ such that, setting 
	$v^{\alpha}(z)=u^{\alpha}(\epsilon_{\alpha}z+z_{\alpha})$,
	for a subsequence, one has:
	\begin{description}
		\item[$(i)$] 
		$\epsilon_{\alpha}\to 0$ and 
		$\epsilon_{\alpha}^{2(\alpha-1)}\to\lambda~$,
	
		\item[$(ii)$] 
		$v^{\alpha}\to\omega$ strongly in $H^{1}_{loc}(\R^{2},\R^{3})$ 
		and uniformly on compact sets of $\R^{2}$,
	
		\item[$(iii)$] 
		$\omega$ is a nonconstant solution to 
		$\Delta\omega=2\lambda H(\omega)\omega_{x}\wedge\omega_{y}$ 
		on $\R^{2}$,
	
		\item[$(iv)$] 
		$\E_{\lambda H}(\omega)\le\lambda\liminf_{\alpha\to 1}
		\E_{H}^{\alpha}(u^{\alpha})$.				
	\end{description}
\end{Proposition}
\eject
\subsection{Proof of Lemma \ref{L:approximating-solutions}}
\noindent
\begin{Lemma}
	\label{L:local-minimum}
	Let $\rho\in(0,(\frac{S_{H}}{2})^{3/2}]$ being $S_{H}$ given by
	{\rm (\ref{eq:isoperimetric-inequality})}.
	Then, for every $u\in H^{1,2\alpha}_{0}$ such that 
	$\|\nabla u\|_{2}\le\rho$ one has 
	$\E_{H}^{\alpha}(u)\ge\frac{1}{2}\D(u)$.
\end{Lemma}

\Proof
Using the inequality $(1+s^{2})^{\alpha}\ge 1+\alpha s^{2}$, one 
infers that $\E_{H}^{\alpha}(u)\ge\E_{H}(u)$ for every $u\in 
H^{1,2\alpha}_{0}$. In addition, by the isoperimetric inequality 
(\ref{eq:isoperimetric-inequality}) one has
$\E_{H}(u)\ge\D(u)-2S_{H}^{-{3}/{2}}\D(u)^{3/2}$.
Therefore $\|\nabla u\|_{2}\le\rho$ implies $\E_{H}^{\alpha}(u)
\ge(1-\sqrt{2}S_{H}^{-{3}/{2}}\rho)\D(u)$ and the thesis follows 
since $\rho\le(\frac{S_{H}}{2})^{3/2}$.
\QED

\begin{Lemma}
	\label{L:uniform-continuity}
	If $u\in H^{1,2\bar\alpha}_{0}$ for some $\bar\alpha>1$, then
	$\E_{H}^{\alpha}(su)\to\E_{H}(su)$ as $\alpha\to 1$, uniformly with 
	respect to $s\in[0,\bar s]$ for every $\bar s>0$.
\end{Lemma}

\Proof
The thesis follows by the estimate
\begin{eqnarray*}
	0&\le&\E_{H}^{\alpha}(su)-\E_{H}(su)=
	\frac{1}{2\alpha}\int_{\disc}\left((1+s^{2}|\nabla u|^{2})^{\alpha}
	-1-\alpha s^{2}|\nabla u|^{2}\right)\\
	&\le&\frac{1}{2\alpha}\int_{\disc}\left(2^{\alpha-1}-1
	+2^{\alpha-1}s^{2\alpha}|\nabla u|^{2\alpha}-\alpha s^{2}
	|\nabla u|^{2}\right)~,
\end{eqnarray*}
and by standard techniques.
\QED

\begin{Lemma}
	\label{L:PS-condition}
	If $\bar M_{H}<+\infty$ then for $\alpha\in(1,\frac{3}{2})$ 
	the functional $\E^{\alpha}_{H}$ 
	satisfies the Palais-Smale condition on $H^{1,2\alpha}_{0}$.
\end{Lemma}

\Proof
First, note that for every $u\in H^{1,2\alpha}_{0}$, using 
(\ref{eq:barMH-definition}) one has
\begin{eqnarray*}
3\E_{H}^{\alpha}(u)-d\E_{H}^{\alpha}(u)u&\ge&
\left(\frac{3}{2\alpha}-1\right)\D^{\alpha}(u)
+2\int_{\disc}(3m_{H}(u)-H(u))u\cdot u_{x}\wedge u_{y}\\
&\ge&\left(\frac{3}{2\alpha}-1\right)\|\nabla u\|_{2\alpha}^{2\alpha}
-\frac{\bar M_{H}}{2}\|\nabla u\|_{2}^{2}~.
\end{eqnarray*}
Hence, 
\begin{equation}
\label{eq:alpha-PS-boundedness}
\left(\frac{3}{2\alpha}-1\right)\|\nabla u\|_{2\alpha}^{2\alpha}
\le\bar M_{H}C_{\alpha}\|\nabla u\|_{2\alpha}^{2}
+\|d\E_{H}^{\alpha}(u)\|
~\|\nabla u\|_{2\alpha}+3\E_{H}^{\alpha}(u)~.
\end{equation}
Now, let $(u^{n})\subset H^{1,2\alpha}_{0}$ be a Palais-Smale 
sequence 
for $\E_{H}^{\alpha}$.
By (\ref{eq:alpha-PS-boundedness}) the sequence $(u^{n})$ is bounded 
in
$H^{1,2\alpha}_{0}$. 
Then, there exists $\bar{u}\in H^{1,2\alpha}_{0}$ such that
(for a subsequence) $u^{n} \to \bar u$ weakly in $H^{1,2\alpha}$ and
uniformly on $\overline\disc$ (by Rellich Theorem).
We need the following auxiliary result (see \cite{[BrCo1]}, for a 
proof):

\begin{Lemma}
	\label{L:V-compactness}
	Let $(u^{n}), (v^{n})\subset H^{1}_{0}\cap L^{\infty}$ be such that 
	$u^{n}\to u$ weakly in $H^{1}$ and $v^{n}\to v$ uniformly. Then
	$$
	\int_{\disc}v^{n}\cdot u^{n}_{x}\wedge u^{n}_{y}\to
	\int_{\disc}v\cdot u_{x}\wedge u_{y}~.
	$$
\end{Lemma}

\noindent
Since for every $h\in H^{1,2\alpha}_{0}$ 
$$
\int_{\disc}(1+|\nabla u^{n}|^{2})^{\alpha-1}\nabla u^{n}\cdot
\nabla h+2\int_{\disc}H(u^{n})h\cdot u^{n}_{x}\wedge 
u^{n}_{y}\to 0
$$
as $n\to+\infty$, thanks to Lemma \ref{L:V-compactness} we obtain 
that $d\E_{H}^{\alpha}(\bar u)=0$.
In particular $0=d\E_{H}^{\alpha}(\bar u)(u^{n}-\bar u) 
=d\D^{\alpha}(u^{n})(u^{n}-\bar u) + o(1)$.
On the other hand, we can use again Lemma \ref{L:V-compactness}
to get $o(1) =d\E^{\alpha}_{H}(u^{n})(u^{n}-\bar u) =
d\D^{\alpha}(u^{n})(u^{n}-\bar u) + o(1)$.
Therefore, $(d\D^{\alpha}(u^{n})-d\D^{\alpha}(\bar u))(u^{n}-\bar 
u)=o(1)$.
Finally we note that $\D^{\alpha}$ is strictly convex on 
$H^{1,2\alpha}_{0}$, and hence $d\D^{\alpha}$ is strictly monotone.
This readily leads to the conclusion.
\QED

In conclusion, we notice that the first part 
of Lemma \ref{L:approximating-solutions} 
is an immediate consequence of Lemmata 
\ref{L:local-minimum} and \ref{L:uniform-continuity}.
The existence of the critical point $u^{\alpha}$ is 
obtained as an application of the mountain pass theorem, 
and by Lemma \ref{L:PS-condition}.

\subsection{Proof of Proposition \ref{P:global-estimates}}
\noindent
In order to show the first estimate, it is useful to introduce, for 
every $\alpha\in(1,\bar\alpha)$, the value
$$
{c}^{\alpha}_{H}=\inf_{u\in H^{1,2\alpha}_{0}\atop u\ne 0}
\sup_{s>0}~\E_{H}^{\alpha}(su)~.
$$

\begin{Lemma}
	\label{L:energy-upper-bound}
	Let $H\in C^{1}(\R^{3})$ satisfy $\bar M_{H}<1$.
	Then $\limsup_{\alpha\to 1}{c}^{\alpha}_{H}\le{c}_{H}$.
\end{Lemma}

\Proof
Fix $\epsilon>0$ and take $u\in C^{\infty}_{c}(\disc,\R^{3})$ such 
that $\sup_{s>0}\E_{H}(su)<{c}_{H}+\epsilon$.
For every $s\ge 0$, using Lemma \ref{L:radial-mp-level}, 
one has
\begin{eqnarray}
	\nonumber
	\E_{H}^{\alpha}(su)&=&\D^{\alpha}(su)-\D(su)+\E_{H}(su)\\
	\label{eq:s-alpha}
	&\le& C_{0}(s^{2\alpha}+1)-C_{1}s^{3}
\end{eqnarray}
with $C_{0},C_{1}>0$ depending just on $u$ (and not on $\alpha$).
Therefore, for $\alpha\in(1,\frac{3}{2})$ there exists
$\bar s_{\alpha}>0$ such that $\E_{H}^{\alpha}(\bar s_{\alpha}u)=
\sup_{s>0}\E_{H}^{\alpha}(su)$. 
{From} (\ref{eq:s-alpha}) it follows that $\bar s_{\alpha}$ is 
uniformly 
bounded.
Then, by Lemma \ref{L:uniform-continuity},
$\lim_{\alpha\to 1}\E_{H}^{\alpha}(\bar s_{\alpha}u)=
\E_{H}(\bar su)$ for some $\bar s>0$.
Hence, $\limsup_{\alpha\to 1}{c}_{H}^{\alpha}\le
{c}_{H}+\epsilon$ and the thesis follows by the 
arbitrariness of $\epsilon>0$.
\QED

Concerning the $H^{1}_{0}$ bounds we have the following result.

\begin{Lemma}
	\label{L:H1-bound}
	Let $H\in C^{1}(\R^{3})$ satisfy $\bar M_{H}<1$.
	If $u\in H^{1,2\alpha}_{0}$ is a nonzero critical point of 
	$\E_{H}^{\alpha}$, then
	$$
	\frac{1}{2}
	\left(\frac{2-\bar M_{H}}{3}\right)^{2}S_{H}^{3}\le
	\int_{\disc}|\nabla u|^{2}\le
	\left(\frac{1}{2\alpha}-\frac{1}{3}-\frac{\bar M_{H}}
	{6}\right)^{-1}
	\E^{\alpha}_{H}(u)~,
	$$
	where $S_{H}$ is given by 
	{\rm (\ref{eq:isoperimetric-inequality})}. 
\end{Lemma}

\Proof
Using (\ref{eq:barMH-definition}) one has
\begin{eqnarray*}
	3\E^{\alpha}_{H}(u)&=&3\D^{\alpha}(u)-d\D^{\alpha}(u)u
	+2\int_{\disc}(3m_{H}(u)-H(u))u\cdot u_{x}\wedge u_{y}\\
	&\ge&
	\left(\frac{3}{2\alpha}-1\right)\int_{\disc}
	\left(1+|\nabla u|^{2}\right)^{\alpha-1}|\nabla u|^{2}+
	\frac{3}{2\alpha}\int_{\disc}\left(
	\left(1+|\nabla u|^{2}\right)^{\alpha-1}-1\right)\\
	& &-\frac{\bar M_{H}}{2}\int_{\disc}|\nabla u|^{2}\\
	&\ge&
	\left(\frac{3}{2\alpha}-1-\frac{\bar M_{H}}{2}\right)
	\int_{\disc}|\nabla u|^{2}~.
\end{eqnarray*}
Moreover, by (\ref{eq:isoperimetric-inequality}) and 
(\ref{eq:barMH-definition}) again, one has
\begin{eqnarray*}
	2\D(u)&\le&
	\int_{\disc}\left(1+|\nabla u|^{2}\right)^{\alpha-1}
	|\nabla u|^{2}\\
	&=&-6\V_{H}(u)+2\int_{\disc}\left(3m_{H}(u)-H(u)\right)
	u\cdot u_{x}\wedge u_{y}\\
	&\le& 6 S_{H}^{-\frac{3}{2}}\D(u)^{\frac{3}{2}}+\bar M_{H}\D(u)~.
\end{eqnarray*}
Since $u\not=0$ one gets the thesis.
\QED

Finally, to show the $L^{\infty}$ bound, $H$ is asked to be constant 
far out and the following estimate holds.

\begin{Lemma}
	\label{L:Linfty-bound}
	Let $H\in C^{1}(\R^{3})$ be such that $H(u)=H_{0}$ for 
	$|u|\ge R_{0}$, where $R_{0}>0$ is given.
	If $u\in H^{1,2\alpha}_{0}$ is a critical point of 
	$\E_{H}^{\alpha}$, then
	$$
	\|u\|_{\infty}\le C|H_{0}|~\|\nabla u\|_{2}^{2}+R_{0}
	$$
	where $C$ is a universal positive constant (independent of
	$\alpha,R_{0},H_{0}$ and $u$).
\end{Lemma}

\Proof
If $u\in H^{1,2\alpha}_{0}$ is a critical point of 
$\E_{H}^{\alpha}$, then $u$ is a weak solution to problem
$$	
\cases{
\div(a_{\alpha}(z)\nabla u)=2H(u)u_{x}\wedge u_{y}
&in $\disc$\cr
u=0&on $\partial\disc$}
$$
where $a_{\alpha}(z)=(1+|\nabla u(z)|^{2})^{\alpha-1}$.
Fix $R>R_{0}$ and let $\Omega_{0}$ be a component of $\{z\in\disc:
|u(z)|>R\}$, if there exists.
Since $u$ is continuous, the set $\Omega_{0}$ is nonempty, bounded, 
open and connected, and $|u|=R$ on $\partial\Omega_{0}$.
Taking $\delta\in(0,R-R_{0})$ one can find a bounded, smooth domain
$\Omega=\Omega_{\delta}$ close to $\Omega_{0}$ such that
$|u(z)|>R_{0}$ for $z\in\Omega$ and $|u(z)|\le R+\delta$ for 
$z\in\partial\Omega$.
Hence $u$ satisfies
\begin{equation}
\label{eq:equation-Ualpha}
\div(a_{\alpha}(z)\nabla u)=2H_{0}
u_{x}\wedge u_{y}\quad{\rm on\ }\Omega~.
\end{equation}
For every $k\in\N$ let $a_{\alpha}^{k}=\min\{a_{\alpha},k\}$ and
let $\varphi^{k}$ be the solution to problem
\begin{equation}
    \label{eq:problem-varphi-alpha-k}
    \cases{
    \div(a_{\alpha}^{k}(z)\nabla\varphi)=g&in $\Omega$\cr
    \varphi=0&on $\partial\Omega$~}
\end{equation}
where $g=2H_{0}u_{x}\wedge u_{y}$.
Since $a_{\alpha}^{k}$ is a continuous bounded function on $\Omega$
and $a_{\alpha}^{k}\ge 1$, by a result of Bethuel and Ghidaglia, 
Theorem 1.3 in \cite{[BeGh]}, there exists a constant $C>0$ such that
\begin{equation}
\label{eq:bound-varphi-alpha-k}
\|\varphi^{k}\|_{\infty}+\|\nabla\varphi^{k}\|_{2}
\le C|H_{0}|~\|\nabla u\|_{2}^{2}
\end{equation}
and $C$ is independent of $k,\alpha,\Omega$ and $u$.
Hence the sequence $(\varphi^{k})$ is bounded in $H^{1}_{0}
(\Omega)$ and thus, there exists $\varphi\in H^{1}_{0}
(\Omega)$ such that, for a subsequence, $\varphi^{k}\to
\varphi$ weakly in $H^{1}_{0}(\Omega)$ and pointwise a.e.
We remark that $a_{\alpha}\in L^{{\alpha\over\alpha-1}}(\Omega)$ 
since $u\in H^{1,2\alpha}_{0}$. 
In particular $a_{\alpha}\in L^{2}(\Omega)$ for $\alpha<2$ and 
$a_{\alpha}^{k}\to a_{\alpha}$ strongly in $L^{2}(\Omega)$.
By (\ref{eq:problem-varphi-alpha-k}) for every
$h\in C^{\infty}_{c}(\Omega)$
$$
\int_{\Omega}a_{\alpha}^{k}(z)\nabla\varphi^{k}
\cdot\nabla h=-\int_{\Omega}gh~.
$$
Hence, by a standard limit procedure, we obtain that for every
$h\in C^{\infty}_{c}(\Omega)$
$$
\int_{\Omega}a_{\alpha}(z)\nabla\varphi
\cdot\nabla h=-\int_{\Omega}gh
$$
that is, $\varphi$ is a weak solution to
\begin{equation}
    \label{eq:problem-varphi-alpha}
    \cases{
    \div(a_{\alpha}(z)\nabla\varphi)=g&in $\Omega$\cr
    \varphi=0&on $\partial\Omega$~.}
\end{equation}
Moreover, by (\ref{eq:bound-varphi-alpha-k}) we also get
\begin{equation}
\label{eq:bound-varphi-alpha}
\|\varphi\|_{\infty}+\|\nabla\varphi\|_{2}
\le C|H_{0}|~\|\nabla u\|_{2}^{2}~.
\end{equation}
Now, we observe that, thanks to (\ref{eq:equation-Ualpha})
and (\ref{eq:problem-varphi-alpha}), the function
$\psi=u-\varphi$ is the solution to problem
\begin{equation}
    \label{eq:problem-psi-alpha}
    \cases{
    \div(a_{\alpha}(z)\nabla\psi)=0&in $\Omega$\cr
    \psi=u&on $\partial\Omega$}
\end{equation}
and it can be characterized as the minimum for the problem
\begin{equation}
    \label{eq:minimization-problem-psi-alpha}
\inf\Big\{\int_{\Omega}a_{\alpha}(z)|\nabla\psi|^{2}:
\psi\in u+H^{1}_{0}(\Omega)\Big\}~.
\end{equation}
Hence $\|\psi\|_{\infty}\le R+\delta$.
Otherwise, if $P$ denotes the projection on the disc
$\disc_{R+\delta}$, that is
$$
P(z)=\cases{
z&if $|z|\le R+\delta$\cr
(R+\delta){z\over|z|}&if $|z|>R+\delta,$}
$$
then $\bar{\psi}=P\circ{\psi}$ will be a
solution to (\ref{eq:minimization-problem-psi-alpha}) and then to
(\ref{eq:problem-psi-alpha}).
In conclusion, using (\ref{eq:bound-varphi-alpha}),
$$
\|u\|_{\infty}\le\|\varphi\|_{\infty}
+\|\psi\|_{\infty}\le
C|H_{0}|~\|\nabla u\|_{2}^{2}+R+\delta~,
$$
and by the arbitrariness of $R>R_{0}$ and $\delta>0$ one gets
the thesis.
\QED

Finally, Proposition \ref{P:global-estimates} 
follows by Lemmata \ref{L:energy-upper-bound}, \ref{L:H1-bound} 
and \ref{L:Linfty-bound}, noting that ${c}_{H}^{\alpha}\ge
\overline{c}_{H}^{\alpha}=\E_{H}^{\alpha}(u^{\alpha})$.
 
\section{Proof of the main theorem}
\noindent
Here we give the proof of Theorem \ref{T:main-result}.
First, as a preliminary result we consider the case in which $H$ is
constant outside a ball (Subsection 4.1). Then, in Subsection 4.2, 
we remove this 
condition, just asking $H$ to be asymptotic to a constant at 
infinity, according to the assumption $\mathbf{(h_{\infty})}$.
\eject
\subsection{Case $H$ constant far out}
\noindent
The results obtained in the previous Sections allow us to deduce 
the existence of an $H$-bubble when the prescribed curvature $H$ 
satisfies $\mathbf{(h_{1})}$ and is constant far out
(this last condition enters in order to guarantee an $L^{\infty}$ 
bound on the approximating solutions).
More precisely, the following result holds.

\begin{Theorem}
	\label{T:preliminary-result}
	Let $H\in C^{1}(\R^{3})$ verify $\mathbf{(h_{1})}$ and the 
	following conditions:
	\begin{description}
		\item[$(i)$] 
		there exists $\bar u\in C^{\infty}_{c}(\R^{2},\R^{3})$ such that 
		$\E_{H}(\bar u)<0$,
	
		\item[$(ii)$] 
		there exists $R_{0}>0$ and $H_{0}\in\R$ such that $H(u)=H_{0}$ 
		as $|u|\ge R_{0}$.
	\end{description}
	Then there exists an $H$-bubble $\omega$  
	such that $\E_{H}(\omega)={c}_{H}$, 
	where ${c}_{H}$ is defined by {\rm (\ref{eq:radial-mp-level})}. 
\end{Theorem}

\begin{Remark}
	\label{R:preliminary-result}
	{\rm Suppose that in the assumption $(ii)$ $H_{0}\ne 0$.
	Then, by Lemma \ref{L:negative-value}, the condition $(i)$ is 
	automatically fulfilled.
	Moreover, in this case problem (\ref{eq:problem}) admits the 
	(trivial) solution $\omega^{0}$ which parametrizes a sphere of 
	radius $|H_{0}|^{-1}$ placed in the region $|u|>R_{0}$.
	However, the additional information on the energy of the $H$-bubble 
	$\omega$ makes meaningful the above result, since if ${c}_{H}<
	\frac{4\pi}{3H_{0}^{2}}$ then $\omega$ is geometrically different
	{from} $\omega^{0}$.}
\end{Remark}

\Proof
{From} the assumptions $(i)$ and $(ii)$, and since $\bar M_{H}<1$,
thanks to Propositions \ref{P:global-estimates} and \ref{P:blow-up-0},
there exists a function $\omega\in X\cap L^{\infty}$ which is a 
$\lambda H$-bubble with $\lambda\in(0,1]$ and 
$\E_{\lambda H}(\omega)\le\lambda{c}_{H}$.
Since $M_{H}<1$, by Proposition \ref{P:minimal-energy} (applied 
with $\lambda H$ instead of 
$H$), $\E_{\lambda H}(\omega)\ge{c}_{\lambda H}$.
Finally, Lemma \ref{L:c-lambda-H} implies $\E_{\lambda H}(\omega)
\ge{c}_{H}$. Then $\lambda=1$ and the Theorem is proved.
\QED

\subsection{General case}
\noindent
Now we want to remove the hypothesis that $H$ is constant far out, by 
requiring just an asymptotic behaviour at infinity as stated by 
$\mathbf{(h_{\infty})}$.
To this aim, we will use the condition $(*)$.
Our argument consists in approximating $H$ with a sequence of 
functions $(H_{n})\subset C^{1}(\R^{3})$ satisfying the hypotheses 
of Theorem \ref{T:preliminary-result} and then passing to the limit 
on 
the sequence $(\omega^{n})$ of the corresponding $H_{n}$-bubbles.
The information on the energies $\E_{H_{n}}(\omega^{n})$ together 
with the condition $(*)$ will permit us to obtain some $L^{\infty}$ 
bound on the sequence $(\omega^{n})$, and then to get the result.

Thus, let us start with the construction of the sequence $(H_{n})$.

\begin{Lemma}
	\label{L:Hn}
	Let $H\in C^{1}(\R^{3})$ satisfying $\mathbf{(h_{\infty})}$ and 
	let $M_{H}$ be defined by {\rm (\ref{eq:MH-definition})}.
	Then there exists a sequence $(H_{n})\subset C^{1}(\R^{3})$ such 
	that:
	\begin{description}
		\item[$(i)$] 
		$H_{n}\to H$ uniformly on $\R^{3}$,
	
		\item[$(ii)$] 
		for every $n\in\N$ there exists $R_{n}>0$ such that 
		$H_{n}(u)=H_{\infty}$ as $|u|\ge R_{n}$,
	
		\item[$(iii)$] 
		$\sup_{u\in\R^{3}}|\nabla H_{n}(u)\cdot u~u|
		:=M_{H_{n}}\le M_{H}$.
	\end{description}
\end{Lemma}

\Proof
It is not restrictive to suppose $H_{\infty}=0$.
Hence, for every $u\in\R^{3}\setminus\{0\}$ one has
$$
{H}(u)=-\int_{1}^{+\infty}\nabla{H}(su)\cdot u~ds~.
$$ 
Let $\chi\in C^{\infty}(\R,[0,1])$ be such that $\chi(r)=1$ as $r\le 
0$, $\chi(r)=0$ as $r\ge 1$ and $|\chi'|\le 2$. 
Given any sequence $r_{n}\to+\infty$ set 
$\chi_{n}(r)=\chi(r-r_{n})$ and
$$
H_{n}(u)=-\int_{1}^{+\infty}\chi_{n}(s|u|)\nabla{H}(su)\cdot u~ds
$$
for every $u\in\R^{3}\setminus\{0\}$.
By continuity, $H_{n}$ is well defined and continuous on $\R^{3}$.
In fact $H_{n}\in C^{1}(\R^{3})$ and for each 
$u\in\R^{3}\setminus\{0\}$
\begin{equation}
	\nabla H_{n}(u)\cdot u=\left.\frac{d}{ds}H_{n}(su)\right|_{s=1}
	=\chi_{n}(|u|)\nabla{H}(u)\cdot u~.
	\label{eq:Hn2}
\end{equation}
Therefore $(iii)$ holds true.
By the definition of $H_{n}$, one has $H_{n}(u)=0$ as $|u|>r_{n}+1$.
Thus $(ii)$ follows, with $R_{n}=r_{n}+1$.
Moreover (\ref{eq:Hn2}) implies $(iii)$.
Now, notice that
\begin{equation}
	\label{eq:Hn3}
	H_{n}(u)=\chi_{n}(|u|){H}(u)+\int_{r_{n}}^{r_{n}+1}\chi'_{n}(t)
	{H}\left(t\frac{u}{|u|}\right)~dt~.
\end{equation}
Setting $\epsilon_{n}=\sup_{|u|\ge r_{n}}|{H}(u)|$, one has 
that 
\begin{eqnarray*}
	\left|\int_{r_{n}}^{r_{n}+1}\chi'_{n}(t)
	{H}\left(t\frac{u}{|u|}\right)~dt\right|&\le&2\epsilon_{n}\\
	\left|(\chi_{n}(|u|)-1){H}(u)\right|&\le&2\epsilon_{n}~.
\end{eqnarray*}
Hence, (\ref{eq:Hn3}) implies that $|H_{n}(u)-{H}(u)|\le 
4\epsilon_{n}$ 
for every $u\in\R^{3}$ and then, since $\epsilon_{n}\to 0$, 
also $(i)$ is 
proved.
\QED

As a further tool, we also need the following result.

\begin{Lemma}
	\label{L:Brezis-Coron}
	Let $(\tilde H_{n})\subset C^{1}(\R^{3})$, $H_{\infty}\in\R$ 
	and $(\tilde\omega_{n})\subset X\cap L^{\infty}$ be such that:
	\begin{description}
		\item[$(i)$]
		$\tilde H_{n}\to H_{\infty}$ uniformly on compact sets,
		
		\item[$(ii)$]
		$\sup_{n}\left(\|\nabla\tilde\omega^{n}\|_{2}+
		\|\tilde\omega^{n}\|_{\infty}\right)<+\infty$,
		
		\item[$(iii)$]
		for every $n\in\N$ the function $\tilde\omega^{n}$ solves 
		$\Delta\tilde\omega^{n}=2\tilde H_{n}(\tilde\omega^{n})
		\tilde\omega^{n}_{x}\wedge\tilde\omega^{n}_{y}$ on $\R^{2}$.
	
	\end{description}
	Then $H_{\infty}\ne 0$ and $\liminf\E_{\tilde H_{n}}
	(\tilde\omega^{n})\ge\frac{4\pi}{3H_{\infty}^{2}}$.
\end{Lemma}

\Proof
{From} the assumption $(ii)$, there exists $\omega\in X\cap 
L^{\infty}$ 
such that, for a subsequence, $\nabla\omega^{n}\to\nabla\omega$ weakly
in $(L^{2}(\R^{2},\R^{3}))^{2}$.
Thanks to the invariance of $H$-systems with respect to dilations,
translations and Kelvin transform, we may also assume that 
$\|\nabla\tilde\omega^{n}\|_{\infty}=|\nabla\tilde\omega^{n}(0)|=1$. 
Then, arguing as in the proof of Proposition \ref{P:blow-up}, 
using the hypotheses $(i)$--$(iii)$, one can 
show that $\omega$ is an $H_{\infty}$-bubble, 
and $\tilde\omega^{n}\to\omega$ strongly in 
$H^{1}_{loc}(\R^{2},\R^{3})$ 
and in $L^{\infty}_{loc}(\R^{2},\R^{3})$. In particular it must be
$H_{\infty}\ne 0$ (there exists no $0$-bubble in $X$).
Furthermore, for every $R>0$, one has
\begin{eqnarray}
	\label{eq:bc1}
	\E_{\tilde H_{n}}(\tilde\omega^{n},D_{R})
	&\to&\E_{H_{\infty}}(\omega,D_{R})\\
	\label{eq:bc2}
	\int_{\partial D_{R}}\tilde\omega^{n}\cdot\frac{\partial
	\tilde\omega^{n}}{\partial\nu}&\to&\int_{\partial D_{R}}
	\omega\cdot\frac{\partial\omega}{\partial\nu}~,
\end{eqnarray}
where, in (\ref{eq:bc1}), we used the notation:
$$
\E_{H}(u,\Omega)=\frac{1}{2}\int_{\Omega}|\nabla u|^{2}+
2\int_{\Omega}m_{H}(u)u\cdot u_{x}\wedge u_{y}~.
$$
Now, fixing $\epsilon>0$, let $R>0$ be such that
$$
|\E_{H_{\infty}}(\omega,\R^{2}\setminus D_{R})|<\epsilon~,\ \
\left(|H_{\infty}|~\|\omega\|_{\infty}+1\right)
\int_{\R^{2}\setminus D_{R}}|\nabla\omega|^{2}<\epsilon~.
$$
Multiplying $\Delta\omega=2H_{\infty}\omega_{x}\wedge\omega_{y}$
by $\omega$ and integrating over $\R^{2}\setminus D_{R}$ we find
$$
\left|\int_{\partial D_{R}}\omega\cdot{\partial\omega \over \partial 
\nu}\right|=\left|\int_{\R^{2}\setminus D_{R}}\left(\omega\cdot
\Delta\omega+|\nabla\omega|^{2}\right)\right|\le\epsilon~.
$$
Then, by (\ref{eq:bc2}), one has that
\begin{equation}
	\left|\int_{\partial D_{R}}\tilde\omega^{n}\cdot\frac{\partial
	\tilde\omega^{n}}{\partial\nu}\right|\le\epsilon+o(1)~.
	\label{eq:bc3}
\end{equation}
Now we multiply $\Delta\tilde\omega^{n} = 
2\tilde H_{n}(\tilde\omega^{n})\tilde\omega^{n}_{x}
\wedge\tilde\omega^{n}_{y}$ by $\tilde\omega^{n}$ 
and we integrate over $\R^{2}\setminus D_{R}$ to get 
\begin{eqnarray}
	\nonumber
	\int_{\partial D_{R}}\tilde\omega^{n}\cdot{\partial\tilde\omega^{n}
	\over\partial\nu}&=&\int_{\R^{2}\setminus D_{R}}
	|\nabla\tilde\omega^{n}|^{2}+2\int_{\R^{2}\setminus D_{R}}
	\tilde H_{n}(\tilde\omega^{n})\tilde\omega^{n}\cdot
	\tilde\omega^{n}_{x}\wedge\tilde\omega^{n}_{y}\\
	\nonumber
	&=&3\E_{\tilde H_{n}}(\tilde\omega^{n},\R^{2}\setminus D_{R})
	-\frac{1}{2}\int_{\R^{2}\setminus D_{R}}
	|\nabla\tilde\omega^{n}|^{2}\\
	\nonumber
	& &+2\int_{\R^{2}\setminus D_{R}}
	(\tilde H_{n}(\tilde\omega^{n})-3m_{\tilde H_{n}}(\tilde\omega^{n}))
	\tilde\omega^{n}\cdot\tilde\omega^{n}_{x}\wedge\tilde\omega^{n}_{y}\\
	\label{eq:bc4}
	&\le&3\E_{\tilde H_{n}}(\tilde\omega^{n},\R^{2}\setminus D_{R})
	-\left(\frac{1}{2}-\mu_{n}\rho\right)
	\int_{\R^{2}\setminus D_{R}}|\nabla\tilde\omega^{n}|^{2}
\end{eqnarray}
where $\rho=\sup_{n}\|\tilde\omega^{n}\|_{\infty}$, and
$\mu_{n}=\sup_{|u|\le\rho}|\tilde H_{n}(u)-3m_{\tilde H_{n}}(u)|$.
Hence (\ref{eq:bc3}) and (\ref{eq:bc4}) imply
$$
\E_{\tilde H_{n}}(\tilde\omega^{n},\R^{2}\setminus D_{R})\ge
-{\epsilon\over 3}+o(1)~,
$$
because, by $(i)$, $\mu_{n}\to 0$.
Finally, we have
\begin{eqnarray*}
\E_{H_{\infty}}(\omega)-\epsilon&\le&
\E_{H_{\infty}}(\omega,D_{R})=\E_{\tilde 
H_{n}}(\tilde\omega^{n},D_{R})
+o(1)\\
&=&\E_{\tilde H_{n}}(\tilde\omega^{n})-\E_{\tilde H_{n}}
(\tilde\omega^{n},\R^{2}\setminus D_{R})+o(1)
\le\E_{\tilde H_{n}}(\tilde\omega^{n}) +\frac{\epsilon}{3} + o(1)~.
\end{eqnarray*}
Hence, by the arbitrariness of $\epsilon>0$, one obtains 
$\liminf\E_{\tilde H_{n}}(\tilde\omega^{n})\ge\E_{H_{\infty}}(\omega)$
and the thesis follows by Remark \ref{R:constant-curvature}.
\QED

\noindent
{\bf Proof of Theorem \ref{T:main-result}}.
Let $(H_{n})\subset C^{1}(\R^{3})$ be the sequence given by Lemma 
\ref{L:Hn}. From Theorem \ref{T:preliminary-result}, for every $n$ 
there exists an $H_{n}$-bubble $\omega^{n}$ such that 
$\E_{H_{n}}(\omega^{n})={c}_{H_{n}}$. 
By Lemma \ref{L:semicontinuity-c-H}, one has
that 
\begin{equation}
	\label{eq:proof1}
	\limsup_{n\to+\infty}\E_{H_{n}}(\omega^{n})\le{c}_{H}~.
\end{equation}
We point out that if we prove that 
$\sup_{n}\|\omega^{n}\|_{\infty}=R<+\infty$, then we have concluded, 
since for $n$ large, $H(u)=H_{n}(u)$ as $|u|\le R$.
To this goal, as a first step, we show that
\begin{equation}
	\label{eq:proof2}
	\|\omega^{n}-\omega^{n}_{\infty}\|_{\infty}\le C_{1}
	\left(1+\int_{\R^{2}}|\nabla\omega^{n}|^{2}\right)
\end{equation}
where $\omega^{n}_{\infty}=\lim_{|z|\to\infty}\omega^{n}(z)$ and
$C_{1}>0$ depends only on $\|H\|_{\infty}$.
This is a consequence of an {\em a~priori} $L^{\infty}$ estimate
proved by Gr\"uter \cite{[Gr]} (see also Theorem 4.8 in 
\cite{[BeRe]}).
More precisely, fixing an arbitrary $\delta>0$, for every $n$ there 
exists $\rho_{n}>0$, depending on $\delta$, such that if 
$|z|\ge\rho_{n}$ then $|\omega^{n}(z)-\omega^{n}_{\infty}|\le\delta$.
Let us set
\begin{eqnarray*}
	\gamma^{n}(z)&=&\omega^{n}(\rho_{n}z)-\omega^{n}_{\infty}
	\ \ {\rm as}\ z\in\partial\disc\\
	u^{n}(z)&=&\omega^{n}(\rho_{n}z)-\omega^{n}_{\infty}
	\ \ {\rm as}\ z\in\disc~.
\end{eqnarray*}
Thus $u^{n}$ is a smooth and conformal solution to
$$
\cases{
\Delta u^{n}=2\tilde H_{n}(u^{n})
u^{n}_{x}\wedge u^{n}_{y}&in $\disc$\cr
u^{n}=\gamma^{n}~&on $\partial\disc~,$ }
$$
where $\tilde H_{n}(u)=H_{n}(u+\omega^{n}_{\infty})$.
Hence, by \cite{[Gr]},
$$
\|u^{n}\|_{L^{\infty}(\disc)}\le
\|\gamma^{n}\|_{L^{\infty}(\partial\disc)}
+C\left(1+\int_{\disc}|\nabla\omega^{n}|^{2}\right)
$$
with $C>0$ depending on $\|\tilde 
H_{n}\|_{\infty}=\|H_{n}\|_{\infty}$.
Since $H_{n}\to H$ uniformly on $\R^{3}$, actually, $C$ is 
independent 
of $n$, but depends only on $\|H\|_{\infty}$.
Then
$$
\|\omega^{n}-\omega^{n}_{\infty}\|_{L^{\infty}(\R^{2})}\le\delta+
\|u^{n}\|_{L^{\infty}(\disc)}\le 2\delta+C
\left(1+\int_{\R^{2}}|\nabla\omega^{n}|^{2}\right)~.
$$
Therefore (\ref{eq:proof2}) holds true.
As a second step, we show that for every $n$
\begin{equation}
	\label{eq:proof3}
	\int_{\R^{2}}|\nabla\omega^{n}|^{2}\le C_{2}
\end{equation}
where $C_{2}>0$ depends only on $H$.
Indeed, by (\ref{eq:MH-definition}), using 
(\ref{eq:omega-condition}), 
one has $(1-M_{H_{n}})\D(\omega^{n})\le 3\E_{H_{n}}(\omega^{n})$.
Since $M_{H_{n}}\le M_{H}$, from (\ref{eq:proof1}) it follows that
$$
\limsup_{n\to+\infty}\int_{\R^{2}}|\nabla\omega^{n}|^{2}\le
\frac{6{c}_{H}}{1-M_{H}}~,
$$
and thus (\ref{eq:proof3}) is proved. Consequently, by 
(\ref{eq:proof2}), one obtains
\begin{equation}
	\label{eq:proof4}
	\|\omega^{n}-\omega^{n}_{\infty}\|_{\infty}\le C_{3}
\end{equation}
with $C_{3}>0$ independent of $n$.
As a last step, let us show that 
$\sup_{n}|\omega^{n}_{\infty}|<+\infty$.
We argue by contradiction, assuming that (for a subsequence) 
$|\omega^{n}_{\infty}|\to+\infty$.
Since $H_{n}\to H$ uniformly on $\R^{3}$, by $\mathbf{(h_{\infty})}$, 
we have that $\tilde H_{n}\to H_{\infty}$ uniformly on compact sets.
Moreover, $\tilde\omega^{n}(z)=\omega^{n}(z)-\omega^{n}_{\infty}$ is 
an $\tilde H_{n}$-bubble and, thanks to (\ref{eq:proof4}) and 
(\ref{eq:proof3}), we can apply Lemma \ref{L:Brezis-Coron}, 
to infer that $H_{\infty}\ne 0$ and
$$
\frac{4\pi}{3H_{\infty}^{2}}\le
\liminf\E_{\tilde H_{n}}(\tilde\omega^{n})=
\liminf\E_{H_{n}}(\omega^{n})~.
$$
Then (\ref{eq:proof1}) implies that 
$\frac{4\pi}{3H_{\infty}^{2}}\le{c}_{H}$, contrary to the 
condition $(*)$.
Therefore, we have that $\sup|\omega^{n}_{\infty}|<+\infty$, 
that, together with (\ref{eq:proof4}), gives the desired estimate.
This concludes the proof.
\QED

We end the work, by making some comments about the case of 
{\em radially symmetric} curvatures.

\begin{Example}
	\label{Ex:radial-curvature}
	{\rm Let $H\in C^{1}(\R^{3})$ be a radial function satisfying 
	$\mathbf{(h_{1})}$ and $\mathbf{(h_{\infty})}$ with $H_{\infty}\ne 
0$.
	Given $\phi\colon\R^{2}\to\S^{2}$ defined by 
	(\ref{eq:stereographic-projection}), and $\rho>0$, 
	the mapping $\rho\phi$ is a solution to (\ref{eq:problem}),
	i.e., it is a radial $H$-bubble, if and only if $\rho|H(\rho)|=1$. 
	In this case the energy of this radial $H$-bubble is 
	$\frac{4\pi}{3H(\rho)^{2}}$.
	Clearly, since $H$ is regular and $H_{\infty}\ne 0$, the 
	equation $\rho|H(\rho)|=1$ always admits positive solutions.
	Now, suppose, in addition that the condition $(*)$ holds true. 
	This happens, for instance, if $H(\rho)>H_{\infty}>0$ for $\rho$ 
	large. Then, there exist $H$-bubbles with minimal energy
	$c_{H}<\frac{4\pi}{3H_{\infty}^{2}}$. Hence, these minimal 
	$H$-bubbles cannot be radial if $|H(\rho)|\le H_{\infty}$ 
	whenever $\rho|H(\rho)|=1$.}
\end{Example}
\nonumsection{Acknowledgments}
Work supported by M.U.R.S.T. progetto di ricerca 
``Metodi Variazionali ed Equazioni Differenziali Nonlineari'' 
(cofin. 2001/02)

\appendix
\section{Appendix}
\subsection{Convergence of approximating solutions in a 
\\ Sacks-Uhlenbeck type setting}
\noindent
In this Appendix we study the behaviour of sequences of 
solutions of approximating problems of the type
$$
\cases{
\div((1+|\nabla u|^{2})^{\alpha-1}\nabla u)=2H(u)u_{x}\wedge u_{y}
&in $\disc$\cr
u=0&on $\partial\disc$}
$$
in the limit as $\alpha\to 1_{+}$.
More precisely, we assume that for every $\alpha\in(1,\bar\alpha)$
a function $u^{\alpha}\in H^{1,2\alpha}_{0}$ is given, 
in such a way that
\begin{eqnarray}
	\label{E:critico}
	& &~~d\E^{\alpha}_{H}(u^{\alpha})=0,\\
	\label{eq:alpha-norms-bound}
	& &\sup_{\alpha\in(1,\bar\alpha)}\left(\|u^{\alpha}\|_{\infty}
	+\|\nabla u^{\alpha}\|_{2}\right)<+\infty~,\\
	\label{eq:alpha-lower-bound}
	& &\inf_{\alpha\in(1,\bar\alpha)}\|\nabla u^{\alpha}\|_{2}>0~.
\end{eqnarray}
The first main result is non-variational and concerns a blow up
analysis of sequences of approximating solutions. 
We point out that this result applies to {\em any} sequence of 
functions satisfying (\ref{E:critico})--(\ref{eq:alpha-lower-bound}).

\begin{Proposition}
	\label{P:blow-up}
	Let $H\in C^{1}(\R^{3})\cap L^{\infty}$ and for every 
	$\alpha\in(1,\bar\alpha)$ 
	let $u^{\alpha}\in H^{1,2\alpha}_{0}$ satisfy 
	{\rm (\ref{E:critico})--(\ref{eq:alpha-lower-bound})}.
	Then, there exist sequences 
	$(\epsilon_{\alpha})\subset(0,+\infty)$, 
	$(z_{\alpha})\subset\overline{\disc}$,
	a number $\lambda\in(0,1]$, and a function 
	$\omega\in X\cap L^{\infty}$ such that, setting 
	$v^{\alpha}(z)=u^{\alpha}(\epsilon_{\alpha}z+z_{\alpha})$,
	for a subsequence, one has:
	\begin{description}
		\item[$(i)$] 
		$\epsilon_{\alpha}\to 0$ and 
		$\epsilon_{\alpha}^{2(\alpha-1)}\to\lambda~$,
	
		\item[$(ii)$] 
		$v^{\alpha}\to\omega$ strongly in $H^{1}_{loc}(\R^{2},\R^{3})$ 
		and uniformly on compact sets of $\R^{2}$,
	
		\item[$(iii)$] 
		$\omega$ is a nonconstant solution to 
		$\Delta\omega=2\lambda H(\omega)\omega_{x}\wedge\omega_{y}$ 
		on $\R^{2}$.		
	\end{description}
\end{Proposition}

Notice that, according to Proposition \ref{P:blow-up}, in the limiting
problem the curvature function is $\lambda H$, with $\lambda\in (0,1]$,
and not necessarily $\lambda=1$.

The second important result of this Appendix is variational and
states a semicontinuity property, under an additional assumption 
on $H$, involving the value $\bar M_{H}$ defined by
{\rm (\ref{eq:barMH-definition})}.

\begin{Proposition}
	\label{P:lambda-estimate}
	Let $H\in C^{1}(\R^{3})$ be such that 
	$\bar M_{H}<1$. For $\alpha\in(1,\bar\alpha)$
	let $u^{\alpha}\in H^{1,2\alpha}_{0}$ satisfy 
	{\rm (\ref{E:critico})--(\ref{eq:alpha-lower-bound})}, and let
	$\lambda\in(0,1]$ and $\omega\in X\cap L^{\infty}$ be given by 
	Proposition \ref{P:blow-up}. Then
	$$
	\E_{\lambda H}(\omega)\le\lambda\liminf_{\alpha\to 1}
	\E_{H}^{\alpha}(u^{\alpha}).
	$$
\end{Proposition}

To prove Proposition \ref{P:blow-up}, first of all we need some local 
estimates on the family $(u^{\alpha})$. This will be developed in 
Subsection A.1. Then the proof of Proposition \ref{P:blow-up} will be 
performed in Subsection A.2. Finally, Proposition 
\ref{P:lambda-estimate} will be proved in Subsection A.3.

\subsection{Local estimates ($\varepsilon$-regularity)}
\label{SS1}
\noindent
Here we study the regularity properties of critical
points for $\E_{H}^{\alpha}$, following the arguments by
Sacks and Uhlenbeck \cite{[SaUh]}.

The first (minor) difference with respect to the framework of
Sacks and Uhlenbeck paper lies in the nonlinear term. In
\cite{[SaUh]} the Euler-Lagrange equation for the
harmonic map problem involves the second fundamental form
of the embedding of the target space $N$ into an Euclidean
space, instead of the curvature term. This is far to lead to
any extra difficulty, since the invariance of the curvature term
with respect to dilations makes computations even easier,
in this case. 

The main difference with \cite{[SaUh]} concerns
the $L^{\infty}$ bound on the maps $u$ under consideration.
In their paper, Sacks and Uhlenbeck deal with maps $u$ whose
target space is a compact Riemannian manifold, and therefore, they
have a natural $L^{\infty}$ bound on all maps $u$. On the
contrary, the target space of our maps $u$ is the noncompact
space $\R^{3}$, and hence we have no natural {\em a priori} 
bound. Therefore we have to ask it as an hypothesis. 

Another difference with respect to the proof of
Sacks and Uhlenbeck is due to the presence of a boundary
in the domain. However, this does not lead any extra difficulty.
One can argue, for example, as in Struwe \cite{[Str2]}, Proposition
2.6.

The first result concerns global regularity for fixed $\alpha>1$,
and it can be obtained as in \cite{[SaUh]}, using Theorem 1.11.$1'$ 
in \cite{[Mo]} and Struwe \cite{[Str2]}, proof of Proposition 2.6, 
for the regularity up to the boundary.

\begin{Lemma}
	\label{L:Morrey-regularity}
	Let $H\in C^{1}(\R^{3})$ and let $u\in H^{1,2\alpha}_{0}$ be 
	a critical point of $\E^{\alpha}_{H}$ for some $\alpha>1$. 
	Then $u$ belongs to $W^{2,q}(D,\R^{3})$ for every 
	$q\in[1,+\infty)$ and solves
	\begin{equation}
		\label{eq:approximate-equation}
		\Delta u=
		-\frac{2(\alpha-1)}{1+|\nabla u|^{2}}
		(\nabla^{2}u,\nabla u)\nabla u
		+\frac{2H(u)}{(1+|\nabla u|^{2})^{\alpha-1}}
		u_{x}\wedge u_{y}\ \ {\it in}\ \disc~.
	\end{equation}
\end{Lemma}

The second result of this Section concerns some local
estimates for the solutions of the approximating problems
({\em $\varepsilon$-regularity}) which are actually the same as
in the celebrated paper \cite{[SaUh]}, and which are stated in the
following Lemma (compare also with Lemma A.1 in \cite{[BeRe]}).
We restrict ourselves to make estimates in the interior of the
disk, thanks to the extension argument by Struwe \cite{[Str2]}.

\begin{Lemma}
	\label{L:main-estimate}
	{\bf (Main Estimate)}
	Let $H\in C^{1}(\R^{3})\cap L^{\infty}$.
	Then there exist $\bar\varepsilon=\bar\varepsilon(\|H\|_{\infty})
	>0$, and for every $p\in(1,+\infty)$ an exponent $\alpha_{p}>1$
	and a constant $C_{p}=C_{p}(\|H\|_{\infty})>0$, 
	such that if $\alpha\in[1,\alpha_{p})$ and $u\in W^{2,p}_{loc}
	(D,\R^{3})$ solves {\rm (\ref{eq:approximate-equation})}, then
	$$
	\|\nabla u\|_{L^{2}(D_{R}(z))}\le\bar\varepsilon
	\ \Rightarrow\
	\|\nabla u\|_{H^{1,p}(D_{R/2}(z))}\le C_{p}R^{\frac{2}{p}-2}
	\|\nabla u\|_{L^{2}(D_{R}(z))}
	$$
	for every disc $\overline{D_{R}(z)}\subset\disc$.
\end{Lemma}

\Proof
Our arguments strictly follow the original proof in
\cite{[SaUh]}.
Let $u$ be a solution to (\ref{eq:approximate-equation}) for some 
$\alpha\ge 1$. 
Fixing $z\in\disc$, for $R\in(0,1-|z|)$ we expand $D_{R}(z)$ to the 
unit 
disc $D$, and we define a map $\omega\colon\disc\to\R^{3}$ by setting
$$
\omega(\zeta)=u(R\zeta+z)-
{\int\!\!\!\!\!\!{-}}_{D_{R}(z)}u
$$
A direct computation shows that $\omega$
is a regular solution in $D$ to the system
\begin{equation}
	\label{E:approximate-omega2}
	\Delta\omega=
	-\frac{2(\alpha-1)}{R^{2}+|\nabla\omega|^{2}}
	(\nabla^{2}\omega,\nabla\omega)\nabla\omega
	+2\frac{H_{R}(\omega)}{(R^{2}+|\nabla\omega|^{2})^{\alpha-1}}
	\omega_{x}\wedge\omega_{y}
\end{equation}
where $H_{R}(\omega)=R^{2(\alpha -1)}H(\omega+
{\int\!\!\!\!\!{-}}_{D_{R}(z)}u)$.
Note also that $R^{-2(\alpha -1)}\|H_{R}\|_{\infty}
\le{\overline H}=:\|H\|_{\infty}$.
Now fix four radii ${1\over 2}=r_{0}<r_{1}< r_{2}<r_{3}=1$ and three 
cut-off functions $\varphi_{i}\in C^{\infty}(\R^{2},[0,1])$
such that $\varphi_{i}\equiv 1$ on $D_{r_{i-1}}$, $\varphi_{i}\equiv
0$ on $\R^{2}\setminus D_{r_{i}}$ ($i=1,2,3$).
Let $K=\max_{i}\left(\|\nabla\varphi_{i}\|_{\infty}+\|\nabla^{2}
\varphi_{i}\|_{\infty}\right)$.
Our aim is to use equation (\ref{E:approximate-omega2}) in
order to obtain some estimate on $\varphi_{i}\omega$.
First, we point out some simple inequalities:
\begin{equation}
	\label{stimo-laplaciano}
	|\Delta(\varphi_{i}\omega)|\le\varphi_{i}|\Delta\omega|+
	2|\nabla\varphi_{i}|~|\nabla\omega|+
	|\Delta\varphi_{i}|~|\omega|~,
\end{equation}
and
\begin{equation}
	\label{stimo1}
	\left|{\varphi_{i}(\nabla^{2}\omega,\nabla\omega)\nabla\omega
	\over R^{2}+|\nabla \omega|^{2}}\right|\le
	|\varphi_{i}\nabla^{2}\omega|\le
	|\nabla^{2}(\varphi_{i}\omega)|
	+2|\nabla\varphi_{i}|~|\nabla\omega|+
	|\nabla^{2}\varphi_{i}|~|\omega|~.
\end{equation}
In order to handle the curvature term in (\ref{E:approximate-omega2})
we observe that
$
2\varphi_{i}(\omega_{x}\wedge\omega_{y})=[(\varphi_{i}\omega)_{x}
\wedge\omega_{y}+\omega_{x}\wedge(\varphi_{i}\omega)_{y}]-
[(\varphi_{i})_{x}(\omega\wedge\omega_{y})+
(\varphi_{i})_{y}(\omega_{x}\wedge\omega)]
$
and hence $|2\varphi_{i}(\omega_{x}\wedge \omega_{y})|
\le 2|\nabla(\varphi_{i}\omega)|~|\nabla \omega|+
|\nabla\varphi_{i}|~|\omega|~|\nabla\omega|$.
Therefore, we can estimate
\begin{equation}
	\label{stimo2}
	\left|{2\varphi_{i}H_{R}(\omega)\omega_{x}\wedge\omega_{y}
	\over(R^{2}+ |\nabla \omega|^{2})^{\alpha-1}}\right|\le
	2{\overline H}|\nabla(\varphi_{i}\omega)|~|\nabla \omega|+
	{\overline H}|\nabla\varphi_{i}|~|\omega|~|\nabla\omega|~.
\end{equation}
Multiplying (\ref{E:approximate-omega2})
by $\varphi_{i}$ and using (\ref{stimo-laplaciano})--(\ref{stimo2})
we obtain
\begin{eqnarray*}
	|\Delta(\varphi_{i}\omega)|&\le&
	2(\alpha-1)|\nabla^{2}(\varphi_{i}\omega)|
	+6K\chi_{i}(|\omega|+|\nabla\omega|)\\
	&\quad&
	+2{\overline H}~|\nabla(\varphi_{i}\omega)|~|\nabla \omega|
	+{\overline H}~|\nabla\varphi_{i}|~|\omega|~|\nabla \omega|
\end{eqnarray*}
where $\chi_{i}$ is the characteristic function of the set
$D_{r_{i}}$.
Thus, for all $p\in(1,+\infty)$ we have
\begin{eqnarray}
	\label{Lp-prima}
	\nonumber
	\|\Delta(\varphi_{i}\omega)\|_{L^{p}(D_{r_{i}})}
	&\!\!\!\le&\!\!\!\!
	2(\alpha-1)\|\varphi_{i}\omega\|_{H^{2,p}(D_{r_{i}})}
	+6K\left(\|\omega\|_{L^{p}(D_{r_{i}})}
	+\|\nabla\omega\|_{L^{p}(D_{r_{i}})}\right)\\
	&\!\!\!+&\!\!\!\!
	2{\overline H}\|~|\nabla(\varphi_{i}\omega)|~|\nabla
	\omega|~\|_{L^{p}(D_{r_{i}})}
	+{\overline H}K\|~|\omega|~|\nabla
	\omega|~\|_{L^{p}(D_{r_{i}})}.
\end{eqnarray}
Since $\omega$ has zero mean value on $D$, we have that for every
$p\in(1,+\infty)$
\begin{equation}
	\label{Lp1}
	\|\omega\|_{L^{p}(D_{r_{i}})}\le C_{p}\|\nabla\omega\|_{L^{2}(D)}~,
\end{equation}
where $C_{p}$ depends only on the Sobolev embedding constant of
$H^{1,2}(D)$ into $L^{p}(D)$ and on the Poincar\'e constant on $D$.
Taking $p\in(1,2]$, we plainly have
\begin{equation}
	\label{eq:trivial}
	\|\nabla\omega\|_{L^{p}(D_{r_{i}})}\le
	2\|\nabla\omega\|_{L^{2}(D)}~.
\end{equation}
Moreover, for $p\in(1,4)$, using H\"older inequality and (\ref{Lp1}),
we can estimate
\begin{eqnarray}
	\label{eq:holder1}
	& &\|~|\nabla(\varphi_{i}\omega)|~|\nabla\omega|~\|_{L^{p}
	(D_{r_{i}})}\le\|\nabla(\varphi_{i}\omega)\|_{L^{4}
	(D_{r_{i}})}\|\nabla\omega\|_{L^{{4p/(4-p)}}(D_{r_{i}})}~,\\
	\label{eq:holder2}
	& &\|~|\omega|~|\nabla\omega|~\|_{L^{p}(D_{r_{i}})}\le
	C_{4}\|\nabla\omega\|_{L^{2}(D)}
	\|\nabla\omega\|_{L^{{4p/(4-p)}}(D_{r_{i}})}~.
\end{eqnarray}
Now we apply the standard regularity theory for linear elliptic
equations.
Denoting by $c(p)$ the norm of the operator $\Delta^{-1}$ as a map
{from} $L^{p}(D_{r_{i}})$ into $W^{2,p}\cap H^{1}_{0}(D_{r_{i}})$,
and using (\ref{Lp-prima})--(\ref{eq:holder2}), we obtain the
following crucial inequality for $p\in(1,2]$
\begin{eqnarray}
	\label{Lp}
	\nonumber
	\beta_{p,\alpha}\|\varphi_{i}\omega\|_{H^{2,p}(D_{r_{i}})}
	&\le&
	{\overline C}_{p}\|\nabla\omega\|_{L^{2}(D)}
	+C_{4}{\overline H}K\|\nabla\omega\|_{L^{2}(D)}
	\|\nabla\omega\|_{L^{{4p/(4-p)}}(D_{r_{i}})}\\
	& &+2{\overline H}~\|\nabla(\varphi_{i}\omega)\|_{L^{4}
	(D_{r_{i}})}\|\nabla\omega\|_{L^{{4p/(4-p)}}(D_{r_{i}})}~.
\end{eqnarray}
where we have set $\beta_{p,\alpha}=c(p)^{-1}-2(\alpha-1)$
and ${\overline C}_{p}=6K(C_{p}+2)$.
First, we use (\ref{Lp}) taking $p=2$ and $i=2$.
{From} (\ref{Lp1}), we have that
$\|\nabla(\varphi_{2}\omega)\|_{L^{4}(D_{r_{2}})}
\le C_{4}K\|\nabla\omega\|_{L^{2}(D)}+
\|\nabla\omega\|_{L^{4}(D_{r_{2}})}$.
Then, we fix $\bar\alpha>1$ such that
$\beta_{2,\bar\alpha}>0$, and we observe that $\bar\alpha$
depends only on the constants in elliptic regularity theory.
Hence, if $\alpha\in[1,\bar\alpha]$, (\ref{Lp}) with $p=2$ 
and $i=2$ yields
\begin{eqnarray}
	\label{eq:p=2}
	\|\omega\|_{H^{2,2}(D_{r_{1}})}
	&\le&\|\varphi_{2}\omega\|_{H^{2,2}(D_{r_{2}})} \\
	\nonumber
	&\le&C_{1}(\overline{H})\left(\|\nabla\omega\|_{L^{2}(D)}
	+\|\nabla\omega\|_{L^{4}(D_{r_{2}})}\|\nabla\omega\|_{L^{2}(D)}
	+\|\nabla\omega\|^{2}_{L^{4}(D_{r_{2}})}\right),
\end{eqnarray}
where $C_{1}(\overline{H})$ depends only on  ${\overline H}$.
Now we show that
$\|\nabla\omega\|_{L^{4}(D_{r_{2}})}$ can be controlled in terms of 
$\|\nabla\omega\|_{L^{2}(D)}$, if $\|\nabla\omega\|_{L^{2}(D)}$ is 
small enough.
To do this, we use again (\ref{Lp}) taking $p={4\over 3}$ and $i=3$.
We point out that the critical Sobolev exponent corresponding to
$p=\frac{4}{3}$ is $p^{*}=4$.
Hence, there exists $S_{4/3}>0$ (independent of the domain) such 
that $S_{4/3}\|\varphi_{3}\omega\|_{H^{1,4}(D)}\le
\|\varphi_{3}\omega\|_{H^{2,4/3}(D)}$.
Therefore, reminding that $r_{3}=1$, (\ref{Lp}), with $p={4\over 3}$ 
and $i=3$, yields
$$
\left(\beta_{{4\over 3},\alpha}-2{\overline H}
S_{4/3}^{-1}\|\nabla\omega\|_{L^{2}(D)}\right)
\|\varphi_{3}\omega\|_{H^{2,4/3}(D)}\le
\overline{C}_{4/3}\|\nabla\omega\|_{L^{2}(D)}
+C_{4}{\overline H}K\|\nabla\omega\|_{L^{2}(D)}^{2}.
$$
Now, take a smaller $\bar\alpha>1$ in order that 
$\beta_{4/3,\bar\alpha}=\bar\beta>0$.
Thus, for every $\alpha\in[1,\bar\alpha]$ we have 
$\beta_{4/3,\alpha}\ge\bar\beta$.
Then, take $\bar\varepsilon>0$ small enough, such that
$\bar\beta - 2S_{4/3}^{-1}{\overline H}\bar\varepsilon>0$.
Notice that $\bar\varepsilon$ depends only on ${\overline H}$.
Therefore, we infer that
\begin{eqnarray*}
	\|\nabla\omega\|_{L^{4}(D_{r_{2}})}&\le&
	\|\nabla(\varphi_{3}\omega)\|_{L^{4}(D)}\le
	\|\varphi_{3}\omega\|_{H^{1,4}(D)}\\
	&\le&S_{4/3}^{-1}\|\varphi_{3}\omega\|_{H^{2,4/3}(D)}\le
	C_{2}(\overline{H})\|\nabla\omega\|_{L^{2}(D)}~,
\end{eqnarray*}
if $\|\nabla \omega\|_{L^{2}(D)}\le\bar\varepsilon$, with 
$C_{2}(\overline{H})$ depending only on $\overline{H}$.
Going back to (\ref{eq:p=2}), we have proved that
$$
\|\omega\|_{H^{2,2}(D_{r_{1}})}\le C_{3}(\overline{H})
\|\nabla\omega\|_{L^{2}(D)}
$$ 
when $\alpha\in[1,\bar\alpha]$,
provided that $\|\nabla\omega\|_{L^{2}(D)}\le\bar\varepsilon$,
being $C_{3}(\overline{H})$ a positive constant depending only on 
$\overline H$.
Hence, by the Sobolev embeddings, for every $q\in[1,+\infty)$ there 
exists a positive constant $C_{4}(q,\overline{H})$, depending also on 
$q$ 
such that 
\begin{equation}
	\label{eq:main-estimate1}
	\|\omega\|_{H^{1,q}(D_{r_{1}})}\le C_{4}(q,\overline{H})
	\|\nabla\omega\|_{L^{2}(D)}
\end{equation}
when $\alpha\in[1,\bar\alpha]$
and $\|\nabla\omega\|_{L^{2}(D)}\le\bar\varepsilon$.
For the last step, we apply (\ref{Lp-prima}) with $i=1$ and we use 
the 
following estimates, obtained with the H\"older inequality and with
(\ref{Lp1}):
\begin{eqnarray*}
	\|~|\nabla(\varphi_{1}\omega)|~|\nabla\omega|~\|_{L^{p}(D_{r_{1}})}
	&\le&K\|~|\omega|~|\nabla\omega|~\|_{L^{p}(D_{r_{1}})}+
	\|\nabla\omega\|^{2}_{L^{2p}(D_{r_{1}})}~,\\
	\|~|\omega|~|\nabla\omega|~\|_{L^{p}(D_{r_{1}})}&\le&C_{2p}
	\|\nabla\omega\|_{L^{2}(D)}\|\nabla\omega\|_{L^{2p}(D_{r_{1}})}~.
\end{eqnarray*}
Then, arguing as for (\ref{Lp}) we get
\begin{eqnarray*}
	\beta_{p,\alpha}\|\varphi_{1}\omega\|_{H^{2,p}(D_{r_{1}})}&\le&
	6K\|\nabla\omega\|_{L^{2}(D)}+6K\|\nabla\omega\|_{L^{p}(D_{r_{1}})}\\
	& &+2\overline{H}\|\nabla\omega\|_{L^{2p}(D_{r_{1}})}
	+3\overline{H}KC_{2p}\|\nabla\omega\|_{L^{2}(D)}
	\|\nabla\omega\|_{L^{2p}(D_{r_{1}})}.
\end{eqnarray*}
Finally, in order to estimate $\|\nabla\omega\|_{L^{p}(D_{r_{1}})}$
and $\|\nabla\omega\|_{L^{2p}(D_{r_{1}})}$, we use 
(\ref{eq:main-estimate1}) with $q=p$ and $q=2p$.
Thus, for fixed $p\in(1,+\infty)$ we can find 
$\alpha_{p}\in(1,\bar\alpha]$ such that for $\alpha\in[1,\alpha_{p}]$ 
one has $\beta_{p,\alpha}\ge\beta_{p,\alpha_{p}}>0$.
Moreover, we can also find a constant $C_{5}(p,\overline{H})>0$ such 
that
for $\alpha\in[1,\alpha_{p}]$, one has
$$
\|\omega\|_{H^{2,p}(D_{1/2})}\le
\|\varphi_{1}\omega\|_{H^{2,p}(D_{r_{1}})}\le
C_{5}(p,\overline{H})\|\nabla\omega\|_{L^{2}(D)}
$$
provided that $\|\nabla\omega\|_{L^{2}(D)}\le\bar\varepsilon$.
To conclude the proof, we just have to remark that
$\|\nabla\omega\|_{L^{2}(D)}=\|\nabla u\|_{L^{2}(D_{R}(z))}$,
and $\|\nabla u\|_{H^{1,p}(D_{R/2}(z))}^{p}=
R^{2-p}\|\nabla\omega\|_{L^{p}(D_{1/2})}^{p}+
R^{2-2p}\|\nabla^{2}\omega\|_{L^{p}(D_{1/2})}^{p}\le
R^{2-2p}\|\omega\|_{H^{2,p}(D_{1/2})}^{p}$, since $R\le 1$.
\QED

\subsection{Passing to the limit (blow up analysis for 
$(u^{\alpha})$)}
\noindent
The first preliminary result concerns the behaviour of the starting 
sequence $(u^{\alpha})$ satisfying 
(\ref{E:critico})--(\ref{eq:alpha-lower-bound}).

\begin{Lemma}
	\label{L:zero-weak-limit}
	$u^{\alpha}\to 0$ weakly in $H^{1}_{0}$ and 
	$\|\nabla u^{\alpha}\|_{\infty}\to+\infty$ as $\alpha\to 1$.
\end{Lemma}

\Proof
Since $(u^{\alpha})$ is bounded in $H^{1}_{0}$ and in
$L^{\infty}$, passing to a 
subsequence, we can assume that $u^{\alpha}\to u$ weakly in 
$H^{1}_{0}$, for some $u\in H^{1}_{0}\cap L^{\infty}$.
Let us prove that $u$ is a weak solution to the Dirichlet problem
\begin{equation}
	\cases{\Delta u=2H(u)u_{x}\wedge u_{y}& in $\disc$\cr
	u=0& on $\partial\disc$}
	\label{eq:dirichlet-problem}
\end{equation}
To this aim, fix an open set $\Omega$ with $\overline{\Omega}
\subset D$.
Arguing as in \cite{[SaUh]}, proof of Proposition 4.3, we can
find a finite set of points $F\subset \Omega$ such that
$\int_{D_{R}(z)}|\nabla u|^{2}\le \bar\varepsilon$ for $z\not\in F$
and $R$ small enough (depending on $z$), where $\bar\varepsilon>0$ 
is given by Lemma \ref{L:main-estimate}.
Then, an application of Lemma \ref{L:main-estimate} gives a uniform
bound for $\|\nabla u^{\alpha}\|_{H^{1,2}(D_{R/2}(z))}$.
Noting also that $(u^{\alpha})$ is bounded in $L^{q}(D)$
for every $q\in[1,+\infty]$, we infer that $(u^{\alpha})$
is bounded in $W^{2,2}(D_{R/2}(z))$, and hence, by Rellich Theorem,
$u^{\alpha} \to u$ strongly in $H^{1}(D_{R/2}(z))$ and in
$L^{\infty}(D_{R/2}(z))$. This is sufficient to conclude that
$u$ is a weak solution to the equation
$\Delta u=2H(u)u_{x}\wedge u_{y}$
in $D_{R/2}(z)$ and hence, since $z$ was arbitrarily
chosen, in $\Omega \setminus F$.
Now we can repeat the proof of Theorem 3.6 in \cite{[SaUh]}. 
Assume for simplicity that $F=\{0\}$. 
Let $\eta \in C^{\infty}(\R,[0,1])$ be
such that $\eta (s) = 0$ for $s\le 1$ and $\eta(s)=1$ for $s\ge 2$,
and set $\eta^{k}(s)=\eta(ks)$. 
Given $h\in C^{\infty}_{c}(\Omega,\R^{3})$ we set
$h^{k}(\zeta)=\eta^{k}(|\zeta|)h(\zeta)$. Notice that
$h^{k}$ can be used as test for $u$ to get
\begin{equation}
	\label{@}
	\int_{\Omega}\nabla u\cdot\nabla h^{k}+
	2 \int_{\Omega}H(u)h^{k}\cdot u_{x}\wedge u_{y}=0~.
\end{equation}
Now, since $h^{k}\to h$ weakly$^{*}$ in $L^{\infty}$,
we get $\int_{\Omega}H(u)h^{k}\cdot u_{x}\wedge u_{y}\to
\int_{\Omega}H(u)h\cdot u_{x}\wedge u_{y}$.
Also, $\int_{\Omega}\nabla u\cdot\nabla h^{k}\to\int_{\Omega}
\nabla u\cdot\nabla h$, since, by H\"older inequality,
$\int_{\Omega}|\nabla u\cdot\nabla \eta^{k}||h|\le C\|\nabla
u\|_{L^{2}(D_{2/k})}=o(1)$ as $k\to+\infty$.
Therefore, (\ref{@}) yields in the limit
$$
\int_{\Omega}\nabla u\cdot\nabla h+
2\int_{\Omega}H(u)h\cdot u_{x}\wedge u_{y}=0
$$
for every test function $h\in C^{\infty}_{c}(\Omega,\R^{3})$,
that is, $u$ solves $\Delta u=2H(u)u_{x}\wedge u_{y}$ in $\Omega$.
Finally, for the arbitrariness of $\Omega$, we conclude that $u$
is a weak solution to problem (\ref{eq:dirichlet-problem}).
Then, by a Heinz regularity result \cite{[He]}, $u$ is smooth, and
a nonexistence result by Wente \cite{[We]}, which holds also in 
case $H$ nonconstant, can be applied, to conclude that $u\equiv 0$.
Thus, we obtain that $u^{\alpha}\to 0$ weakly in $H^{1}_{0}$ and
strongly in $H^{1}_{loc}(\disc\setminus N)$ where $N$ is a countable 
set of $\disc$. In particular $\nabla u^{\alpha}\to 0$ pointwise a.e. 
in $\disc$. Therefore, as a last step, we observe that if it were
$\liminf_{\alpha\to 1}\|\nabla u^{\alpha}\|_{\infty}<+\infty$, then
$\liminf_{\alpha\to 1}\|\nabla u^{\alpha}\|_{2}=0$, contrary to
(\ref{eq:alpha-lower-bound}). Hence, it must be 
$\|\nabla u^{\alpha}\|_{\infty}\to +\infty$ as $\alpha\to 1$.
\QED

\noindent
{\bf Proof of Proposition \ref{P:blow-up}.}
\noindent
For every $\alpha\in(1,\bar\alpha)$ set
$\epsilon_{\alpha}=\|\nabla u^{\alpha}\|_{\infty}^{-1}$, 
let $z_{\alpha}\in\overline\disc$ be such that $|\nabla 
u^{\alpha}(z_{\alpha})|=\epsilon_{\alpha}^{-1}$ and define
\begin{equation}
	v^{\alpha}(z)=u^{\alpha}(\epsilon_{\alpha}z+z_{\alpha})~.
	\label{eq:v-alpha-definition}
\end{equation}
Note that $v^{\alpha}\in H^{1}_{0}(D_{\alpha},\R^{3})$ where
$D_{\alpha}=D_{\epsilon_{\alpha}^{-1}}\big(-\frac{z_{\alpha}}
{\epsilon_{\alpha}}\big)$.
Moreover the following facts 
hold:
\begin{eqnarray}
	\label{eq:normalized-L-infty-norm}
	& &\|v^{\alpha}\|_{\infty}=\|u^{\alpha}\|_{\infty}\\
	\label{eq:normalized-gradient-norm}
	& &\|\nabla v^{\alpha}\|_{2}=\|\nabla u^{\alpha}\|_{2}\\
	\label{eq:max-gradient-norm}
	& &|\nabla v^{\alpha}(0)|=\|\nabla v^{\alpha}\|_{\infty}=1~.
\end{eqnarray}
Furthermore, $v^{\alpha}\in W^{2,q}_{loc}(D_{\alpha},\R^{3})$ for 
every $q\in[1,+\infty)$ and solves the system
\begin{equation}
	\label{eq:pb-v-alpha}
	\Delta v^{\alpha}=-\frac{2(\alpha-1)}
	{\epsilon_{\alpha}^{2}+|\nabla v^{\alpha}|^{2}}
	(\nabla^{2}v^{\alpha},\nabla v^{\alpha})\nabla v^{\alpha}
	+\frac{2\epsilon_{\alpha}^{2(\alpha-1)}H(v^{\alpha})}
	{(\epsilon_{\alpha}^{2}+|\nabla v^{\alpha}|^{2})^{\alpha-1}}
	v^{\alpha}_{x}\wedge v^{\alpha}_{y}
	\ \ {\rm in}\ D_{\alpha}~.
\end{equation}
Since $\epsilon_{\alpha}\to 0$ as $\alpha\to 1$, one has that
$0<\epsilon_{\alpha}^{2(\alpha-1)}<1$, and then, for a subsequence,
$\epsilon_{\alpha}^{2(\alpha-1)}\to\lambda$ for some 
$\lambda\in[0,1]$. Moreover, setting 
$\rho_{\alpha}=\epsilon_{\alpha}^{-1}
{\rm dist}(z_{\alpha},\partial\disc)$, we may also assume that there 
exists $\lim_{\alpha\to 1}\rho_{\alpha}\in[0,+\infty]$.
Let $\Omega_{\infty}$ be the union of all compact sets in $\R^{2}$ 
contained in $D_{\alpha}$ as $\alpha\to 1$. Note that
$\Omega_{\infty}$ is a half-plane if 
$\rho_{\alpha}\to\ell\in[0,+\infty)$, 
while $\Omega_{\infty}=\R^{2}$ if $\rho_{\alpha}\to+\infty$.
{From} (\ref{eq:alpha-norms-bound}), 
(\ref{eq:normalized-L-infty-norm})
and (\ref{eq:normalized-gradient-norm}) it follows that there exists
$\omega\in X\cap L^{\infty}$ such that, for a subsequence, $\nabla 
v^{\alpha}\to\nabla\omega$ weakly in $(L^{2}(\R^{2},\R^{3}))^{2}$.
Moreover, by (\ref{eq:max-gradient-norm}) one has that $v^{\alpha}\to
\omega$ strongly in $L^{\infty}_{loc}(\R^{2},\R^{3})$.
Let $\bar\varepsilon>0$ be given by Lemma \ref{L:main-estimate}.
Take an arbitrary compact set $K$ in $\Omega_{\infty}$ and set
$R_{K}={\rm dist}(K,\partial\Omega_{\infty})$. Then, let $R\in
(0,\min\{1,R_{K},\frac{\bar\varepsilon}{\sqrt{\pi}}\big\})$.
Hence, there exists $\alpha_{K}>1$ such that $K\subset
D_{\alpha}$ for $\alpha\in(1,\alpha_{K})$ and, consequently,
for every $z\in K$, one has $\overline{D_{R}(z)}\subset D_{\alpha}$ 
and $\|\nabla v^{\alpha}\|_{2}\le\bar\varepsilon$. 
Because of the definition (\ref{eq:v-alpha-definition}) of 
$v^{\alpha}$, 
one can apply Lemma \ref{L:main-estimate}, in order to conclude that 
$\|\nabla v^{\alpha}\|_{H^{1,p}(D_{R/2}(z))}$ is uniformly bounded 
with 
respect to $\alpha\in(1,\alpha_{K})$, for every $p>1$. 
Using (\ref{eq:normalized-L-infty-norm}) and 
(\ref{eq:alpha-norms-bound}), 
we infer that $(v^{\alpha})$ is bounded in $H^{2,p}(D_{R/2}(z))$. 
Therefore we can conclude that $\omega\in H^{2,p}(D_{R/2}(z))$, 
$v^{\alpha}\to\omega$ strongly in $H^{1}(D_{R/2}(z))$, and $\nabla 
v^{\alpha}\to\nabla\omega$ pointwise everywhere in $D_{R/2}(z)$.
Since $z$ is an arbitrary point in $K$ and $K$ is any compact set 
in $\Omega_{\infty}$, a standard diagonal argument yields that
$\omega\in H^{2,p}_{loc}(\Omega_{\infty})$ for every $p<+\infty$,
$v^{\alpha}\to\omega$ strongly in $H^{1}_{loc}(\Omega_{\infty})$, 
and $\nabla v^{\alpha}\to\nabla\omega$ pointwise everywhere in 
$\R^{2}$.
In particular, by (\ref{eq:max-gradient-norm}), $|\nabla\omega(0)|=
\|\nabla\omega\|_{\infty}=1$, and thus $\omega$ is nonconstant.
Now we test (\ref{eq:pb-v-alpha}) on an arbitrary function $h\in 
C^{\infty}_{c}(D_{R/2}(z),\R^{3})$ and we pass to the limit as 
$\alpha\to 1$. First, we have 
\begin{equation}
	\int_{\R^{2}}\Delta v^{\alpha}\cdot h\to\int_{\R^{2}}\nabla\omega
	\cdot\nabla h~,
	\label{eq:lambda-H-solution-1}
\end{equation}
because of the weak convergence $\nabla v^{\alpha}\to\nabla\omega$.
Secondly, using the estimate
$$
\left|\int_{\R^{2}}\frac{(\nabla^{2}v^{\alpha},\nabla v^{\alpha})
\nabla v^{\alpha}\cdot h}{\epsilon_{\alpha}^{2}+|\nabla 
v^{\alpha}|^{2}}
\right|\le\int_{\R^{2}}|\nabla^{2}v^{\alpha}|~|h|\le
\|\nabla^{2}v^{\alpha}\|_{L^{p}(D_{R/2}(z))}\|h\|_{L^{p'}}
$$
and the fact that $v^{\alpha}$ is uniformly bounded in 
$H^{2,p}(D_{R/2}(z))$ as $\alpha\in(1,\alpha_{K})$, we obtain that
\begin{equation}
	2(\alpha-1)\int_{\R^{2}}\frac{(\nabla^{2}v^{\alpha},
	\nabla v^{\alpha})\nabla v^{\alpha}\cdot h}
	{\epsilon_{\alpha}^{2}+|\nabla v^{\alpha}|^{2}}\to 0,
	\label{eq:lambda-H-solution-2}
\end{equation}
as $\alpha\to 1$.
Lastly, setting 
$$
w^{\alpha}=\epsilon_{\alpha}^{2(\alpha-1)}\left(\frac{1}
{(\epsilon_{\alpha}^{2}+|\nabla v^{\alpha}|^{2})^{\alpha-1}}
-1\right)v^{\alpha}_{x}\wedge v^{\alpha}_{y}
$$
one has
$$
\int_{\R^{2}}\frac{\epsilon_{\alpha}^{2(\alpha-1)}H(v^{\alpha})}
{(\epsilon_{\alpha}^{2}+|\nabla v^{\alpha}|^{2})^{\alpha-1}}
h\cdot v^{\alpha}_{x}\wedge v^{\alpha}_{y}
=\int_{\R^{2}}H(v^{\alpha})h\cdot w^{\alpha}
+\epsilon_{\alpha}^{2(\alpha-1)}\int_{\R^{2}}H(v^{\alpha})
h\cdot v^{\alpha}_{x}\wedge v^{\alpha}_{y}~.
$$
Since $\epsilon_{\alpha}^{2(\alpha-1)}\to\lambda$,
$H(v^{\alpha})\to H(\omega)$ uniformly on 
$\overline{D_{R/2}(z)}$ and $\nabla v^{\alpha}\to\nabla\omega$ 
pointwise in $D_{R/2}(z)$, by (\ref{eq:max-gradient-norm}),
on one hand we infer that 
$$
\epsilon_{\alpha}^{2(\alpha-1)}\int_{\R^{2}}H(v^{\alpha})
h\cdot v^{\alpha}_{x}\wedge v^{\alpha}_{y}\to
\lambda\int_{\R^{2}}H(\omega)h\cdot\omega_{x}\wedge\omega_{y}~.
$$
On the other hand, since $\epsilon_{\alpha}\in(0,1)$, we observe that
$$
|w^{\alpha}|\le(1+\epsilon_{\alpha}^{2(\alpha-1)})|v^{\alpha}_{x}|~
|v^{\alpha}_{y}|\le|\nabla v^{\alpha}|^{2}\le 1
$$
and $w^{\alpha}(\zeta)\to 0$ for every $\zeta\in D_{R/2}(z)$.
Indeed, if $\nabla\omega(\zeta)=0$ then $|w^{\alpha}|\le
|\nabla v^{\alpha}|^{2}\to 0$, while if $\nabla\omega(\zeta)\ne 0$ 
then ${(\epsilon_{\alpha}^{2}+|\nabla 
v^{\alpha}(\zeta)|^{2})^{\alpha-1}}
\to 0$.
In conclusion, by the dominated convergence Theorem, we obtain
that $\int_{\R^{2}}H(v^{\alpha})h\cdot w^{\alpha}\to 0$ and then
\begin{equation}
	\int_{\R^{2}}\frac{2\epsilon_{\alpha}^{2(\alpha-1)}H(v^{\alpha})}
	{(\epsilon_{\alpha}^{2}+|\nabla v^{\alpha}|^{2})^{\alpha-1}}
	h\cdot v^{\alpha}_{x}\wedge v^{\alpha}_{y}\to 2\lambda\int_{\R^{2}}
	H(\omega)h\cdot\omega_{x}\wedge\omega_{y}~,
	\label{eq:lambda-H-solution-3}
\end{equation}
as $\alpha\to 1$.
Then (\ref{eq:pb-v-alpha})--(\ref{eq:lambda-H-solution-3}) imply that
$$	
\int_{\R^{2}}\nabla\omega\cdot\nabla h+2\lambda\int_{\R^{2}}
H(\omega)h\cdot\omega_{x}\wedge\omega_{y}=0
$$
for every $h\in C^{\infty}_{c}(D_{R/2}(z),\R^{3})$, for every 
$z\in K$ and for every compact set $K$ in $\Omega_{\infty}$, that is,
$\omega$ solves $\Delta\omega=2\lambda H(\omega)\omega_{x}\wedge
\omega_{y}$ in $\Omega_{\infty}$.
Suppose that $\Omega_{\infty}$ is a half-plane. Since $v^{\alpha}=0$ 
on $\partial D_{\alpha}$, one has that $\omega=0$ on $\partial
\Omega_{\infty}$.
Moreover, since a half-plane is conformally equivalent to a disc,
$\omega$ gives arise to a nonconstant solution to the Dirichlet 
problem
\begin{equation}
	\cases{\Delta u=2\lambda H(u)u_{x}\wedge u_{y}& in $\disc$\cr
	u=0& on $\partial\disc.$}
	\label{eq:lambda-Dirichlet-problem}
\end{equation}
As already noted in the proof of Lemma \ref{L:zero-weak-limit}, the 
only solution to (\ref{eq:lambda-Dirichlet-problem}) is $u\equiv 0$,
and this gives a contradiction, since $\omega$ is nonconstant.
Hence, it must be $\Omega_{\infty}=\R^{2}$, that is, $\omega$ is a
$\lambda H$-bubble. Finally, we observe that $\lambda>0$, since the 
only bounded solutions to $\Delta u=0$ on $\R^{2}$ with 
$\D(u)<+\infty$ are the constant functions, and we already know that
$\omega$ is nonconstant.
This concludes the proof.
\QED

\subsection{Proof of Proposition \ref{P:lambda-estimate}}
\noindent
For every domain $\Omega$ in $\R^{2}$, $\alpha\in(1,\bar\alpha)$,
and $\lambda\in(0,1]$, set 
\begin{eqnarray*}
	\tilde\E_{H}^{\alpha}(v^{\alpha},\Omega)&=&
	\frac{1}{2\alpha}\int_{\Omega}\left((\epsilon_{\alpha}^{2}
	+|\nabla v^{\alpha}|^{2})^{\alpha}-\epsilon_{\alpha}^{2\alpha}
	\right)+2\epsilon_{\alpha}^{2(\alpha-1)}\int_{\Omega}
	m_{H}(v^{\alpha})v^{\alpha}\cdot v^{\alpha}_{x}\wedge 
v^{\alpha}_{y}\\
	\E_{\lambda H}(\omega,\Omega)&=&
	\frac{1}{2}\int_{\Omega}|\nabla\omega|^{2}
	+2\lambda\int_{\Omega}m_{H}(\omega)
	\omega\cdot\omega_{x}\wedge\omega_{y}~.
\end{eqnarray*}
Notice that $\tilde\E_{H}^{\alpha}(v^{\alpha},D_{\alpha})=
\epsilon_{\alpha}^{2(\alpha-1)}\E_{H}^{\alpha}(u^{\alpha})$ and
$v^{\alpha}$ solves the system
\begin{equation}
	\div\left(
	(\epsilon_{\alpha}^{2}+|\nabla v^{\alpha}|^{2})^{\alpha-1}
	\nabla v^{\alpha}\right)=2\epsilon_{\alpha}^{2(\alpha-1)}
	H(v^{\alpha})v^{\alpha}_{x}\wedge v^{\alpha}_{y}~.
	\label{eq:v-alpha}
\end{equation}
Now, multiplying (\ref{eq:v-alpha}) by $v^{\alpha}$,
we obtain
\begin{eqnarray}
	\nonumber
	\div\left(
	(\epsilon_{\alpha}^{2}+|\nabla v^{\alpha}|^{2})^{\alpha-1}
	\nabla v^{\alpha}\cdot v^{\alpha}\right)&=&
	(\epsilon_{\alpha}^{2}+|\nabla v^{\alpha}|^{2})^{\alpha-1}
	|\nabla v^{\alpha}|^{2}\\
	\label{eq:da-integrare}
	& &+2\epsilon_{\alpha}^{2(\alpha-1)}H(v^{\alpha})
	v^{\alpha}\cdot v^{\alpha}_{x}\wedge v^{\alpha}_{y}~.
\end{eqnarray}
Integrating (\ref{eq:da-integrare}) on a domain $\Omega$ and using 
the divergence theorem we infer that
\begin{eqnarray}
	\nonumber
	\int_{\partial\Omega}
	(\epsilon_{\alpha}^{2}+|\nabla v^{\alpha}|^{2})^{\alpha-1}
	v^{\alpha}\cdot\frac{\partial v^{\alpha}}{\partial\nu}&=&
	\int_{\Omega}(\epsilon_{\alpha}^{2}+|\nabla 
v^{\alpha}|^{2})^{\alpha-1}
	|\nabla v^{\alpha}|^{2}\\
	\label{eq:integrata}
	& &+2\epsilon_{\alpha}^{2(\alpha-1)}\int_{\Omega}H(v^{\alpha})
	v^{\alpha}\cdot v^{\alpha}_{x}\wedge v^{\alpha}_{y}~.
\end{eqnarray}
Using (\ref{eq:barMH-definition}) and the definition of 
$\tilde\E_{H}^{\alpha}(v^{\alpha},\Omega)$ we can estimate
\begin{eqnarray}
	\nonumber
	2\epsilon_{\alpha}^{2(\alpha-1)}\int_{\Omega}H(v^{\alpha})
	v^{\alpha}\cdot v^{\alpha}_{x}\wedge v^{\alpha}_{y}
	&\le&\epsilon_{\alpha}^{2(\alpha-1)}\frac{\bar M_{H}}{2}
	\int_{\Omega}|\nabla v^{\alpha}|^{2}
	+3\tilde\E_{H}^{\alpha}(v^{\alpha},\Omega)\\
	\label{eq:solita-stima}
	& &-\frac{3}{2\alpha}\int_{\Omega}\left((\epsilon_{\alpha}^{2}
	+|\nabla v^{\alpha}|^{2})^{\alpha}-\epsilon_{\alpha}^{2\alpha}
	\right).
\end{eqnarray}
Hence, setting
\begin{eqnarray*}
	I_{\alpha}(\partial\Omega)&=&\frac{1}{3}\int_{\partial\Omega}
	(\epsilon_{\alpha}^{2}+|\nabla v^{\alpha}|^{2})^{\alpha-1}
	v^{\alpha}\cdot\frac{\partial v^{\alpha}}{\partial\nu}\\
	I_{\alpha}(\Omega)&=&
	\frac{1}{2\alpha}\int_{\Omega}\left((\epsilon_{\alpha}^{2}
	+|\nabla v^{\alpha}|^{2})^{\alpha}-\epsilon_{\alpha}^{2\alpha}
	\right)-\frac{1}{3}\int_{\Omega}(\epsilon_{\alpha}^{2}+
	|\nabla v^{\alpha}|^{2})^{\alpha-1}|\nabla v^{\alpha}|^{2}\\
	& &-\epsilon_{\alpha}^{2(\alpha-1)}\frac{\bar M_{H}}{6}
	\int_{\Omega}|\nabla v^{\alpha}|^{2},
\end{eqnarray*}
by (\ref{eq:solita-stima}) the equation (\ref{eq:integrata}) becomes
\begin{equation}
	\label{eq:sintetica}
	\tilde\E_{H}^{\alpha}(v^{\alpha},\Omega)\ge
	I_{\alpha}(\partial\Omega)+I_{\alpha}(\Omega)~.
\end{equation}
With algebraic computations, one has
\begin{eqnarray*}
	I_{\alpha}(\Omega)&\ge&
	\left(\frac{1}{2\alpha}-\frac{1}{3}\right)
	\int_{\Omega}\left((\epsilon_{\alpha}^{2}
	+|\nabla v^{\alpha}|^{2})^{\alpha}-\epsilon_{\alpha}^{2\alpha}
	\right)-\epsilon_{\alpha}^{2(\alpha-1)}\frac{\bar M_{H}}{6}
	\int_{\Omega}|\nabla v^{\alpha}|^{2}\\
	&\ge&\epsilon_{\alpha}^{2(\alpha-1)}
	\left(\frac{1}{2}-\frac{\alpha}{3}-\frac{\bar M_{H}}{6}\right)
	\int_{\Omega}|\nabla v^{\alpha}|^{2}~.
\end{eqnarray*}
Then, since $\bar M_{H}<1$, one obtains that $I_{\alpha}(\Omega)\ge 0$ 
for $\alpha>1$ sufficiently close to 1, whatever $\Omega$ is.
Hence, (\ref{eq:sintetica}) reduces to
\begin{equation}
	\label{eq:naylon}
	\tilde\E_{H}^{\alpha}(v^{\alpha},\Omega)\ge
	I_{\alpha}(\partial\Omega)~.
\end{equation}
Now we take $\Omega=\R^{2}\setminus D_{R}$.
First, we observe that, since $v^{\alpha}\to\omega$ strongly in 
$H^{1}_{loc}(\R^{2},\R^{3})$ and uniformly on compact sets, and
$\epsilon_{\alpha}^{2(\alpha-1)}\to\lambda$, it holds that
\begin{eqnarray*}
	\lim_{\alpha\to 1}\tilde\E_{H}^{\alpha}(v^{\alpha},D_{R})
	&=&\E_{\lambda H}(\omega,D_{R})\\
	\limsup_{\alpha\to 1}\left|I_{\alpha}(\partial D_{R})\right|
	&\le&\frac{1}{3}\left|\int_{\partial D_{R}}\omega\cdot
	\frac{\partial\omega}{\partial\nu}\right|
\end{eqnarray*}
for every $R>0$.
Then, by (\ref{eq:naylon}), we obtain
\begin{eqnarray*}
	\E_{\lambda H}(\omega,D_{R})&=&
	\tilde\E_{H}^{\alpha}(v^{\alpha})
	-\tilde\E_{H}^{\alpha}(v^{\alpha},\R^{2}\setminus D_{R})+o(1)\\
	&\le&\epsilon_{\alpha}^{2(\alpha-1)}\E_{H}^{\alpha}(u^{\alpha})
	+\frac{1}{3}\left|\int_{\partial D_{R}}\omega\cdot
	\frac{\partial\omega}{\partial\nu}\right|+o(1)
\end{eqnarray*}
where $o(1)\to 0$ as $\alpha\to 1$, for every $R>0$.
Hence
\begin{equation}
	\lambda\liminf_{\alpha\to 1}\E_{H}^{\alpha}(u^{\alpha})\ge
	\E_{\lambda H}(\omega,D_{R})-\frac{1}{3}\left|\int_{\partial D_{R}}
	\omega\cdot\frac{\partial\omega}{\partial\nu}\right|
	\label{eq:terital}
\end{equation}
for every $R>0$.
Finally, notice that
\begin{eqnarray*}
	\left|\int_{\partial D_{R}}\omega\cdot
	\frac{\partial\omega}{\partial\nu}\right|&=&
	\left|\int_{\R^{2}\setminus D_{R}}\left(\omega\cdot
	\Delta\omega+|\nabla\omega|^{2}\right)\right|\\
	&=&\left|\int_{\R^{2}\setminus D_{R}}\left(2\lambda H(\omega)
	\omega\cdot\omega_{x}\wedge\omega_{y}+|\nabla\omega|^{2}\right)
	\right|\\
	&\le&\left(\lambda\|H\|_{\infty}\|\omega\|_{\infty}+1\right)
	\int_{\R^{2}\setminus D_{R}}|\nabla\omega|^{2}.
\end{eqnarray*}
Then, passing to the limit as $R\to+\infty$, from (\ref{eq:terital}) 
the thesis follows.
\QED

\nonumsection{References}

\end{document}